\documentclass[preprint,3p,12pt]{elsarticle}
\usepackage{mathrsfs}
\usepackage{amsmath}
\usepackage{stmaryrd}
\usepackage{bbding}
\usepackage{dcolumn}
\usepackage{graphicx}
\usepackage{amsfonts}
\usepackage{amssymb}
\usepackage{psfrag}
\usepackage{wrapfig}
\usepackage{subfigure}
\usepackage{makeidx}
\usepackage{bm}
\usepackage{epsf}
\usepackage{epsfig}
\usepackage{setspace}
\usepackage{graphicx}
\usepackage{epstopdf}
\usepackage{psfrag}
\usepackage{subfigure}
\usepackage{color}
\usepackage{amsthm}

\begin{document}
\title{Finite Difference Implementation of a High-order Space-Time Coupled Compact Gas-Kinetic Scheme}

\author[HKUST1]{Fengxiang Zhao}
\ead{fzhaoac@connect.ust.hk}

\author[IAPCM]{Yibing Chen}
\ead{chen\_yibing@iapcm.ac.cn}

\author[HKUST1,HKUST2]{Kun Xu\corref{cor}}
\ead{makxu@ust.hk}

\address[HKUST1]{Department of Mathematics, Hong Kong University of Science and Technology, Clear Water Bay, Kowloon, Hong Kong, China}
\address[HKUST2]{Shenzhen Research Institute, Hong Kong University of Science and Technology, Shenzhen, China}
\address[IAPCM]{Institute of Applied Physics and Computational Mathematics and National Key laboratory of Computational physics, and Center for Applied Physics and Technology, Peking University, Beijing, China}
\cortext[cor]{Corresponding author}

\begin{abstract}
This study presents a high-order compact finite difference gas-kinetic scheme (FD-CGKS) that introduces a novel spatial discretization strategy for the efficient implementation of space-time coupled high-order schemes on structured grids. A conservative nonlinear compact discretization is achieved by formulating numerical fluxes from physical fluxes at both nodal and interfacial locations. To simplify the multidimensional spatial reconstruction required for the GKS flux evaluation, we propose a dual-grid approach that updates conservative variables on both a primary grid and an identical dual grid, offset by half the mesh spacing. By leveraging the time-accurate interface solutions from the gas-kinetic evolution model, the scheme explicitly updates averaged spatial derivatives between virtual interfaces, naturally enabling compact high-order reconstruction. Furthermore, a nonlinear GENO method is incorporated to capture flow discontinuities with high resolution and robustness, effectively suppressing spurious oscillations. The proposed framework, which also offers new perspectives for designing schemes based on space-time decoupled Riemann solvers, is systematically validated. Comprehensive benchmark computations of inviscid and viscous flows demonstrate the scheme's high accuracy in resolving a wide spectrum of flow features, from smooth multiscale structures to strong shock discontinuities.
\end{abstract}

\begin{keyword}
Compact GKS, Compact finite difference, Dual-grid approach, Space-time coupled method
\end{keyword}

\maketitle

\section{Introduction}

A fundamental challenge in the numerical simulation of compressible flows within computational fluid dynamics (CFD) is simultaneously achieving high resolution for multiscale turbulent structures and robust, non-oscillatory treatment of shock discontinuities. This dual requirement has long motivated the development of high-accuracy, shock-capturing numerical schemes. Over the past several decades, various advanced methods have been established to address this challenge. Notable examples include weighted essentially non-oscillatory (WENO) schemes \cite{liu-WENO,jiang-WENO}, compact finite difference schemes \cite{lele}, dispersion-relation-preserving (DRP) schemes \cite{DRP}, discontinuous Galerkin (DG) methods \cite{cockburn1,cockburn2}, the flux reconstruction/correction procedure via reconstruction (FR/CPR) framework \cite{huynh2007fr,CPR_wang}, and the compact gas-kinetic scheme (GKS) \cite{zhao2019compact,zhao2023direct}.

Building upon the standard GKS, the compact GKS has been systematically developed to achieve up to eighth-order accuracy on structured meshes \cite{zhao2019compact,zhao2020acoustic}. This scheme has demonstrated favorable spectral-like resolution properties alongside strong robustness in benchmark problems involving two-dimensional strong shocks and aeroacoustic phenomena. Subsequently, a unified and rigorous numerical framework for the compact GKS was established through the consistent modeling of the temporal evolution of flow variables and their high-order numerical fluxes at cell interfaces \cite{zhao2023direct}.
More recently, an effective implementation strategy was proposed to simplify the reconstruction procedure and further enhance the resolution on structured meshes~\cite{yang2026-CGKS}, thereby extending the scale-resolving capability of the compact GKS to high-Reynolds-number turbulent flows.
In parallel, the formulation and application of the compact GKS on unstructured meshes have also been actively investigated \cite{zhao_compact-tri,zhao2023AIA}.
Despite these advances, all existing compact GKS formulations are based on the finite volume method. For computations on structured meshes, the finite difference (FD) method generally offers superior computational efficiency, owing to its simpler high-order spatial discretization. Nevertheless, the FD formulation of the compact GKS remains largely unexplored, representing a notable gap in the development of this class of methods.

To address this limitation and enable highly efficient computations on structured meshes, we propose a novel FD discretization framework for the space-time coupled GKS. Specifically, a compact, high-order FD scheme is developed by constructing nonlinear numerical fluxes at virtual interfaces adjacent to a grid node, utilizing physical fluxes from both neighboring nodes and these virtual interfaces. This specific formulation of numerical fluxes inherently guarantees the strict conservation property of the scheme. Furthermore, the construction seamlessly integrates the generalized essentially non-oscillatory (GENO) method \cite{zhao2025generalized} to robustly capture shocks without inducing spurious oscillations.
This study demonstrates that the proposed numerical flux is analogous to the flux function defined in the conventional FD-WENO schemes \cite{jiang-WENO}, where the physical flux is interpreted as the sliding average of the numerical flux function. However, compared to the FD-WENO schemes \cite{jiang-WENO}, the present approach introduces two key distinctions: (i) the numerical flux reconstruction achieves compactness by simultaneously incorporating values from both virtual interfaces and grid nodes, and (ii) it directly employs physical fluxes in the reconstruction, thereby circumventing the need for upwind flux splitting procedures. 
Finally, it is worth noting that a compact FD discretization using both neighboring interfacial physical fluxes and nodal flux derivatives has been proposed within the weighted compact nonlinear scheme (WCNS) framework \cite{deng2000-WCNS}. The present FD scheme nevertheless remains fundamentally distinct from that approach, as it relies solely on nodal flux values rather than on their derivatives.

A further key contribution of this work is the development of a dual-grid approach, which facilitates efficient, high-order, and compact reconstruction on multidimensional grids.
This strategy allows the space-time coupled compact GKS to derive the multidimensional initial states required for genuinely multidimensional flux evaluations by relying exclusively on one-dimensional procedures.
Specifically, variables from both the primary and dual grids are synthesized to determine these initial states via a one-dimensional compact reconstruction. This process leverages both nodal values and the averaged spatial gradients defined across adjacent virtual interfaces. Consequently, the approach preserves the stencil compactness required for both spatial reconstruction and flux differencing. Ultimately, the proposed methodology establishes an efficient and highly generalizable pathway for integrating space-time coupled schemes into FD frameworks.

This paper is organized as follows.
The high-order FD-CGKS formulation for compressible flow simulations will be introduced in Section \ref{sec:1d-scheme}.
Section \ref{sec:reconstruction} presents the compact spatial reconstructions.
Section \ref{sec:2d-scheme} presents the two-dimensional FD-CGKS.
Section \ref{sec:examples} provides validation test cases, and Section \ref{sec:conclusion} concludes the paper.
The Appendix presents the detailed formulations of the GKS and the GENO reconstruction.

\section{One-dimensional formulation}
\label{sec:1d-scheme}

This section presents the one-dimensional FD formulation of the compact GKS. First, the linear discretization is introduced, followed by the development of a nonlinear discretization designed to capture flow discontinuities.

\subsection{Conservative compact finite difference discretization}

This study develops a conservative FD method within a space-time coupled numerical framework for solving compressible flow equations. Although the current work is based on the space-time coupled compact GKS, the proposed methodology is highly generalizable and not limited to this specific framework. The one-dimensional conservation law for compressible flows is expressed as
\begin{equation}\label{eq-conservation-law}
\mathbf{W}_t + \mathbf{F}_x = 0,
\end{equation}
where \(\mathbf{W} = (\rho, \rho U, \rho E)^\mathrm{T}\) is the vector of conservative variables, and \(\mathbf{F}\) represents the corresponding flux vector. In this study, both inviscid and viscous fluxes are taken into account simultaneously.

\begin{figure}[!htb]
\centering
\includegraphics[width=0.85\textwidth]{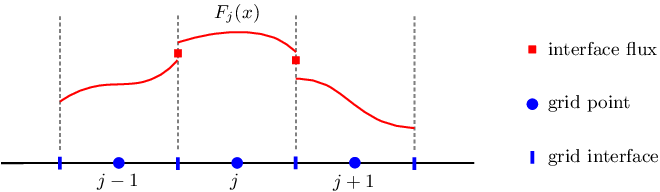}
\caption{\label{1-differ-stencil} Schematic of the stencil for the compact high-order FD discretization of the flux derivative $F_{j,x}$ in FD-CGKS. Both nodal and interfacial fluxes are utilized to achieve the compact discretization.}
\end{figure}

The semi-discrete conservative FD scheme for Eq. \eqref{eq-conservation-law} is formulated as
\begin{equation}\label{fd-conservative-form}
\begin{split}
\frac{\partial \mathbf{W}_j}{\partial t}&=-\frac{1}{\Delta x}\big( \overline{\mathbf{F}}_{j+1/2} -\overline{\mathbf{F}}_{j-1/2} \big),
\end{split}
\end{equation}
where \(\overline{\mathbf{F}}_{j+1/2}\) denotes the numerical flux at the virtual grid interface \(j+1/2\) (hereafter simply referred to as the grid interface). To achieve high-order spatial accuracy, the numerical flux \(\overline{\mathbf{F}}_{j+1/2}\) is reconstructed via a high-order combination of adjacent physical fluxes.
For the overall discretization of the spatial flux derivative, the compact stencil illustrated in Fig. \ref{1-differ-stencil} is employed.
Specifically, this stencil incorporates physical fluxes evaluated at both the grid interfaces, \(\widehat{\mathbf{F}}_{j+1/2}(t)\), and the grid nodes, \(\widehat{\mathbf{F}}_j(t)\).
In the present work, these physical fluxes are evaluated using the flux solver of the GKS, which is based on the gas distribution function.
Although these physical fluxes are inherently time-dependent, the temporal argument \(t\) is omitted in subsequent expressions for notational brevity.

The numerical flux function \(\overline{\mathbf{F}}(x)\) is obtained via a high-order reconstruction on a compact stencil, utilizing the physical fluxes evaluated at adjacent nodes and interfaces. This reconstruction requires that the spatial averages of \(\overline{\mathbf{F}}(x)\) over the local intervals exactly reproduce the corresponding physical fluxes; namely,
\begin{equation}\label{fd-conservative-flux-recons}
\begin{split}
&\frac{1}{\Delta x}\int_{x_{k-1/2}}^{x_{k+1/2}}\overline{\mathbf{F}}(x)\,\mathrm{d}x
=\widehat{\mathbf{F}}_{k},\qquad k=j,\,j+1,\\[2pt]
&\frac{1}{\Delta x}\int_{x_{i-1}}^{x_{i}}\overline{\mathbf{F}}(x)\,\mathrm{d}x
=\widehat{\mathbf{F}}_{i-1/2},\qquad i=j,\,j+1,\,j+2.
\end{split}
\end{equation}
Although the construction of the numerical flux shares similar underlying concepts with conventional FD-WENO schemes~\cite{jiang-WENO}, the present approach distinguishes itself in two key aspects: the high-order reconstruction is achieved on a compact stencil, and the numerical flux is formulated directly from the physical fluxes.
Solving Eq.~\eqref{fd-conservative-flux-recons} and evaluating the reconstructed function at the interface \(x_{j+1/2}\) yields the interface value of the numerical flux,
\begin{equation}\label{fd-conservative-flux}
\overline{\mathbf{F}}_{j+1/2}
=\frac{1}{30}\big(\widehat{\mathbf{F}}_{j-1/2}+46\,\widehat{\mathbf{F}}_{j+1/2}
+\widehat{\mathbf{F}}_{j+3/2}-9(\widehat{\mathbf{F}}_{j}+\widehat{\mathbf{F}}_{j+1})\big).
\end{equation}
In general, this provides a fifth-order accurate approximation of \(\overline{\mathbf{F}}(x_{j+1/2})\). However, owing to the symmetry of the stencil with respect to the interface \(x_{j+1/2}\), the leading \(\mathcal{O}(\Delta x^{5})\) error term cancels out, thereby elevating the formal accuracy to sixth order. A Taylor expansion around \(x_{j+1/2}\) shows that the truncation error contains only even powers of \(\Delta x\),
\[
\overline{\mathbf{F}}_{j+1/2}
=\overline{\mathbf{F}}(x_{j+1/2})
+\frac{1}{8960}\,\Delta x^{6}\,\overline{\mathbf{F}}^{(6)}(x_{j+1/2})
+\mathcal{O}(\Delta x^{8}).
\]

At time \(t^n\), the initial values \(\mathbf{W}^n(x)\) are assumed to be piecewise smooth, with discontinuities confined exclusively to the cell interfaces. This permits the evaluation of physical fluxes at both grid nodes and interfaces: nodal fluxes are computed directly from the smooth local values, whereas interface fluxes are resolved from the discontinuous states to produce a spatially continuous flux. As previously mentioned, all such fluxes are derived from the GKS gas distribution function (\ref{appendix:A}).

Furthermore, wider stencils can be employed to achieve higher-order FD discretizations, for example, by incorporating the values at $j-2$ and $j+2$. The implementation of this approach within the proposed scheme is straightforward and will not be further explored in the present study.

\subsection{Nonlinear numerical flux}

In scenarios where discontinuities such as strong shocks appear within the stencil of Fig. \ref{1-differ-stencil}, large gradients can induce numerical oscillations despite the physical fluxes being theoretically continuous. To suppress such oscillations and robustly capture discontinuities, the reconstruction of the interface numerical fluxes incorporates a nonlinear methodology.

In this study, we propose a nonlinear approach that directly combines a high-order numerical flux with a lower-order flux. Consequently, the nonlinear high-order numerical flux, \(\widetilde{\mathbf{F}}_{j+1/2}\), is formulated as follows:
\begin{equation}\label{fd-conservative-flux-1}
\widetilde{\mathbf{F}}_{j+1/2} = (1-\chi_{j+1/2})\widehat{\mathbf{F}}_{j+1/2} + \chi_{j+1/2} \overline{\mathbf{F}}_{j+1/2}.
\end{equation}
Here, the nonlinear weight \(\omega_{j+1/2}\) is determined by the path function in the GENO nonlinear spatial reconstruction. This weight approaches \(1\) in smooth regions and tends to \(0\) near discontinuities. The detailed formulation of \(\omega_{j+1/2}\) is detailed in the subsequent section on high-order spatial reconstruction. As observed in Eq.~\eqref{fd-conservative-flux-1}, the numerical flux reduces to the purely physical flux in regions containing discontinuities. The scheme with the nonlinear numerical flux given by Eq.~\eqref{fd-conservative-flux-1} is referred to as FD-CGKS(6-2).

Furthermore, an alternative nonlinear numerical flux that reverts to a fourth-order form at discontinuities is given by:
\begin{equation}\label{fd-conservative-flux-2}
\begin{split}
\widetilde{\mathbf{F}}_{j+1/2}= (1-\chi_{j+1/2})\big(8\widehat{\mathbf{F}}_{j+1/2}-\widehat{\mathbf{F}}_{j}-\widehat{\mathbf{F}}_{j+1} \big) + \chi_{j+1/2} \overline{\mathbf{F}}_{j+1/2}.
\end{split}
\end{equation}
The scheme with the nonlinear numerical flux given by Eq. \eqref{fd-conservative-flux-2} is denoted as FD-CGKS(6-4).
Numerical experiments demonstrate that this higher-order formulation achieves a superior balance between numerical accuracy and robustness.

\subsection{Evaluation of the Averaged Gradient in FD-CGKS}

\begin{figure}[!htb]
\centering
\includegraphics[width=0.8\textwidth]{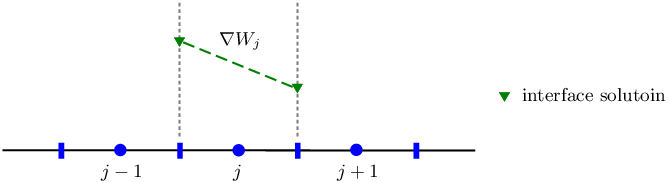}
\caption{\label{1-schematic-gradient} Schematic of the explicit evaluation of averaged gradients over the intervals bounded by adjacent grid interfaces in the FD-CGKS. These averaged gradients, alongside the nodal values, are employed to achieve compact spatial reconstruction.}
\end{figure}

In the FD formulation of the compact GKS, the flow variables at the grid nodes are updated alongside their averaged derivatives over the corresponding intervals bounded by adjacent grid interfaces. By leveraging the evolution solution of the gas distribution function in the GKS, the time-dependent flow variables at the grid interfaces are modeled. Consequently, the averaged gradient is explicitly obtained via Gauss's theorem as follows:
\begin{equation}\label{define-averaged-gradient}
\nabla \mathbf{W}_j(t) = \frac{1}{\Delta x} \big( \mathbf{W}_{j+1/2}(t) - \mathbf{W}_{j-1/2}(t) \big).
\end{equation}
A schematic illustration of the averaged gradient evaluation is presented in Fig. \ref{1-schematic-gradient}. The methodology for evaluating the evolution solution at the grid interfaces within the compact GKS framework was originally proposed in \cite{zhao2023direct}, and the specific computational formulas are detailed in \ref{appendix:A}. Consequently, by simultaneously yielding both the physical fluxes \(\mathbf{F}\) and the macroscopic flow variables \(\mathbf{W}\) at the grid interfaces, the spatio-temporally coupled GKS solver serves as the fundamental building block for constructing the compact FD scheme in this study.

\subsection{High-order time integration}

Within the present finite difference framework, the nonlinear two-stage fourth-order (S2O4) time integration scheme is employed to advance the semi-discrete conservation laws (Eq. \eqref{fd-conservative-form}) and update the nodal variables \(\mathbf{W}_j\). For a given time step \([t^n, t^{n+1}]\), the temporal evolution is given by:
\begin{equation}\label{eq-nodal-update}
\begin{aligned}
\mathbf{W}_j^{n+1/2} &= \mathbf{W}_j^n + \frac{1}{2}\Delta t \mathcal{L}_j^n + \frac{1}{8}\Delta t^2 \frac{\partial}{\partial t}\mathcal{L}_j^n, \\
\mathbf{W}_j^{n+1} &= \mathbf{W}_j^n + \Delta t \mathcal{L}_j^n + \frac{1}{2}\Delta t^2 \frac{\partial}{\partial t}\mathcal{L}_j^n + \frac{1}{3}\Delta t^2 \frac{\partial}{\partial t}(\widetilde{\mathcal{L}}_j^{n+1/2}-\widetilde{\mathcal{L}}_j^n),
\end{aligned}
\end{equation}
where \(\mathcal{L}_j\) denotes the spatial finite difference operator evaluated at node \(j\), defined as \(\mathcal{L}_j(t)=-\big( \overline{\mathbf{F}}_{j+1/2}(t) -\overline{\mathbf{F}}_{j-1/2}(t) \big)/\Delta x\). Additionally, \(\widetilde{\mathcal{L}}_j^n\) represents the associated nonlinear operator. For a comprehensive description of these operators, readers are referred to \cite{zhao2023direct}.

Concurrently, the averaged gradients are advanced to high-order temporal accuracy by employing the same intermediate time stage \(t^{n+1/2}\), yielding:
\begin{equation}\label{eq-interface-update}
\begin{aligned}
\nabla \mathbf{W}_{j}^{n+1/2} &= \nabla \mathbf{W}_{j}^n + \frac{1}{2}\Delta t \, \frac{\partial}{\partial t} \nabla\mathbf{W}_{j}^n, \\
\nabla \mathbf{W}_{j}^{n+1}   &= \nabla \mathbf{W}_{j}^n + \Delta t \, \frac{\partial}{\partial t} \nabla\mathbf{W}_{j}^{n+1/2}.
\end{aligned}
\end{equation}
Here, the averaged gradient \(\nabla \mathbf{W}_{j}^n\) and its temporal derivative \(\frac{\partial}{\partial t} \nabla \mathbf{W}_{j}^n\) are evaluated through Eq.~\eqref{define-averaged-gradient}. The time-dependent evolution solutions of the flow variables at the grid interfaces, which are required for this evaluation, are obtained directly by taking moments of the gas distribution function provided by the GKS. The detailed formulations are presented in \ref{appendix:A}.

\section{Compact spatial reconstruction}
\label{sec:reconstruction}

This section presents a one-dimensional, high-order compact reconstruction designed to determine the flow states at both the grid interfaces and the grid nodes. These reconstructed states are then utilized within the gas-kinetic model to evaluate the numerical fluxes and macroscopic flow variables. Specifically, this compact spatial reconstruction relies on the nodal values \(\mathbf{W}_{j}\) and the averaged gradients \(\mathbf{W}_{j,x}\) defined between adjacent grid interfaces.

\subsection{High-order compact reconstruction}
The compact stencil for the flow state reconstruction at the grid interface \(j+1/2\) is defined as
\[
\overline{S}_{j+1/2} = \big\{ I_{j-1}, I_{j}, I_{j+1}, I_{j+2} \big\},
\]
which encompasses a total of eight known data values. To enhance the numerical stability of the scheme, upwind-biased flow information is employed to reconstruct the left and right states, respectively:
\begin{equation*}
\begin{aligned}
\overline{S}^{\mathrm{data},\mathrm{l}}_{j+1/2} &= \big\{ \mathbf{W}_{j+i}, \nabla \mathbf{W}_{j+k} \big\}, \quad i=-1,0,1,2, \quad k=-1,0,1,\\
\overline{S}^{\mathrm{data},\mathrm{r}}_{j+1/2} &= \big\{ \mathbf{W}_{j+i}, \nabla \mathbf{W}_{j+k} \big\}, \quad i=-1,0,1,2, \quad k=0,1,2.
\end{aligned}
\end{equation*}

Let \(Q\) denote an arbitrary characteristic component of \(\mathbf{W}\). The reconstruction procedure is performed component-wise, achieving seventh-order spatial accuracy. The resulting reconstructions at the interface are given by
\begin{equation}\label{7th-compact-recons}
\begin{aligned}
P^{6,\mathrm{l}}(x_{j+1/2}) &= \frac{1}{64}\big[(7Q_{j-1}+35Q_{j}+21Q_{j+1}+Q_{j+2}) + (Q_{j-1,x}+22Q_{j,x}-7Q_{j+1,x})\Delta x\big],\\
P^{6,\mathrm{l}}_x(x_{j+1/2}) &= -\frac{1}{30\Delta x}\big[(Q_{j-1}+45Q_{j}-45Q_{j+1}-Q_{j+2}) + 9(Q_{j,x}+Q_{j+1,x})\Delta x\big],\\
P^{6,\mathrm{r}}(x_{j+1/2}) &= \frac{1}{64}\big[(Q_{j-1}+21Q_{j}+35Q_{j+1}+7Q_{j+2}) + (7Q_{j,x}-22Q_{j+1,x}-Q_{j+2,x})\Delta x\big],\\
P^{6,\mathrm{r}}_x(x_{j+1/2}) &= -\frac{1}{30\Delta x}\big[(Q_{j-1}+45Q_{j}-45Q_{j+1}-Q_{j+2}) + 9(Q_{j,x}+Q_{j+1,x})\Delta x\big].
\end{aligned}
\end{equation}

Similarly, the compact stencil for the flow state reconstruction at node \(j\) is given by
\[
\overline{S}_{j} = \big\{ I_{j-1}, I_{j}, I_{j+1} \big\},
\]
which provides up to sixth-order spatial accuracy. At the nodal points, only the spatial derivatives require evaluation via reconstruction. Consequently, the resulting sixth-order reconstruction at node \(j\) is obtained by
\begin{equation}\label{6th-compact-recons}
\begin{aligned}
P^{5}(x_{j}) &= Q_j,\\
P^{5}_x(x_{j}) &= \frac{1}{30\Delta x}\big[ 9(Q_{j-1} - Q_{j+1}) + \Delta x(Q_{j-1,x} + 46Q_{j,x} + Q_{j+1,x}) \big].
\end{aligned}
\end{equation}

\subsection{Compact GENO reconstruction}

The GENO method is a nonlinear high-order reconstruction scheme that combines a linear high-order reconstruction with a reliable and physically consistent lower-order reconstruction. This method enables adaptive high-order accuracy for smooth solutions while reverting to the lower-order reconstruction when necessary to preserve numerical robustness. For a comprehensive analysis of the GENO method, readers are referred to \cite{zhao2025generalized}.

The GENO formulation is given by
\begin{equation}\label{GENO-nonlinear}
\begin{aligned}
R(\mathbf{x}) &= \chi P^H(\mathbf{x}) + (1-\chi) P^L(\mathbf{x}),\\
\chi &= \frac{\tanh(C \alpha)}{\tanh(C)},
\end{aligned}
\end{equation}
where \(\chi\) is the path function, and the parameter \(C\) is fixed at \(20\). The ultimate smoothness indicator \(\alpha\), which quantifies the smoothness of the linear high-order reconstruction, is defined as
\begin{equation}\label{GENO-IS-linear}
\begin{aligned}
\alpha &= \frac{2\alpha^H}{\alpha^H+\alpha^L},\\
\alpha^H &= 1 + \left(\frac{IS^{\tau}}{IS^H+\epsilon}\right)^r, \quad \alpha^L = 1 + \left(\frac{IS^{\tau}}{IS^L+\epsilon}\right)^r,
\end{aligned}
\end{equation}
where the exponent \(r\) is set to \(2\), and the small positive constant \(\epsilon\) is taken as \(10^{-12}\) to prevent division by zero. The terms \(IS^H\) and \(IS^L\) are derived from the smoothness indicators associated with the lower-order polynomials \(P^L(\mathbf{x})\) \cite{zhao2025generalized}. \(IS^{\tau}\) serves as a smoothness metric associated with the higher-order derivative terms of the reconstruction over the large stencil.

For the reconstruction at the grid interface \(j+1/2\), the selection of sub-stencils is analogous to that in the CGKS within the finite volume framework. Specifically, the sub-stencils for the left-state reconstruction at the interface are defined as follows:
\begin{equation*}
s_1 = \{\mathbf{W}_{j-1}, \mathbf{W}_{j}, \mathbf{W}_{x,j-1}\}, \quad
s_2 = \{\mathbf{W}_{j-1}, \mathbf{W}_{j}, \mathbf{W}_{j+1}\}, \quad
s_3 = \{\mathbf{W}_{j}, \mathbf{W}_{j+1}, \mathbf{W}_{j+2}\}.
\end{equation*}
The sub-stencils for the right-state reconstruction can be readily obtained by symmetry and are thus omitted here. For the reconstruction at the nodal point \(j\), the following sub-stencils are adopted:
\begin{equation*}
s_1 = \{\mathbf{W}_{j-1}, \mathbf{W}_{j}, \mathbf{W}_{j+1}\}, \quad
s_2 = \{\mathbf{W}_{j-1}, \mathbf{W}_{j}, \mathbf{W}_{x,j-1}\}, \quad
s_3 = \{\mathbf{W}_{j}, \mathbf{W}_{j+1}, \mathbf{W}_{x,j+1}\}.
\end{equation*}
Based on these sub-stencils, the corresponding quadratic polynomials $q^{l,r}_k(x)$ can be easily derived by enforcing exact interpolation of the data values on each stencil; the explicit expressions are provided in \ref{appendix:C}.
The lower-order reconstruction \(P^L\) is constructed as a convex combination of these polynomials:
\begin{equation*}
\begin{aligned}
&P^L(x) = \sum_{k=1}^{3} \omega_k q_k(x), \\
&\omega_k = \frac{\widetilde{\omega}_k}{\sum_{m=1}^{3}\widetilde{\omega}_m}, \quad \widetilde{\omega}_k = \frac{d_k}{(IS_k+\epsilon)^2},
\end{aligned}
\end{equation*}
where the ideal weights are set to \(d_1=8\) and \(d_2=d_3=1\), and the small positive constant is \(\epsilon=10^{-6}\).

Moreover, the smoothness indicators for these reconstruction polynomials are evaluated using their standard definitions \cite{jiang-WENO}. For instance, the smoothness indicator for the left-state reconstruction at the interface \(j+1/2\) is given by
\begin{equation}\label{IS-highorder-1}
IS_k = \sum_{k=1}^{2} \Delta x^{2k-1} \int_{x_{j-1/2}}^{x_{j+1/2}} \left( \frac{\mathrm{d}^{k} q^l_k(x)}{\mathrm{d}x^k} \right)^2 \mathrm{d}x.
\end{equation}
Finally, the resulting \(IS^L\), \(IS^H\), and \(IS^{\tau}\) for the GENO method, which apply to both the interface and nodal reconstructions, are formulated as
\begin{equation}\label{GENO-1d-IS-HL}
\begin{aligned}
\mathit{IS}^L &= \min_{k=1,2,3} \mathit{IS}_k, \quad
\mathit{IS}^H = \max_{k=1,2,3} \mathit{IS}_k, \\
\mathit{IS}^{\tau} &= \left| \frac{\mathit{IS}_2 + \mathit{IS}_3}{2} - \mathit{IS}_1 \right|.
\end{aligned}
\end{equation}

Additionally, the path function at the grid interface also serves to determine the nonlinear weight in the nonlinear numerical flux of Eq.~\eqref{fd-conservative-flux-1} and Eq.~\eqref{fd-conservative-flux-2}, which is given by
\[
\chi_{j+1/2} = \min\{\chi^{l,Q}_{j+1/2}, \chi^{r,Q}_{j+1/2}\},
\]
where \(Q = \rho\) when reconstructing conservative variables, whereas \(Q\) represents the characteristic variable associated with the characteristic speed \(U \pm c\) when reconstructing characteristic variables.

\section{Multidimensional FD formulation}
\label{sec:2d-scheme}

This section presents the multidimensional formulation of the FD-CGKS, focusing primarily on the two-dimensional compact reconstruction, which represents a core innovation of this study. While the extension to three dimensions is straightforward, the present study focuses exclusively on the two-dimensional case. The semi-discrete FD formulation of the hyperbolic conservation laws on a two-dimensional Cartesian grid is given by
\begin{equation}\label{eq-conservation-law-2d}
\mathbf{W}_t = -\frac{1}{\Delta x}\big( \overline{\mathbf{F}}_{i+1/2,j} - \overline{\mathbf{F}}_{i-1/2,j} \big) - \frac{1}{\Delta y}\big( \overline{\mathbf{G}}_{i,j+1/2} - \overline{\mathbf{G}}_{i,j-1/2} \big),
\end{equation}
where \( \mathbf{W} = (\rho, \rho U, \rho V, \rho E)^T \), and \( \overline{\mathbf{F}} \) and \( \overline{\mathbf{G}} \) represent the corresponding flux vectors in the \( x \)- and \( y \)-directions, respectively. Both inviscid and viscous fluxes are taken into account simultaneously.

\begin{figure}[!htb]
\centering
\includegraphics[width=0.60\textwidth]{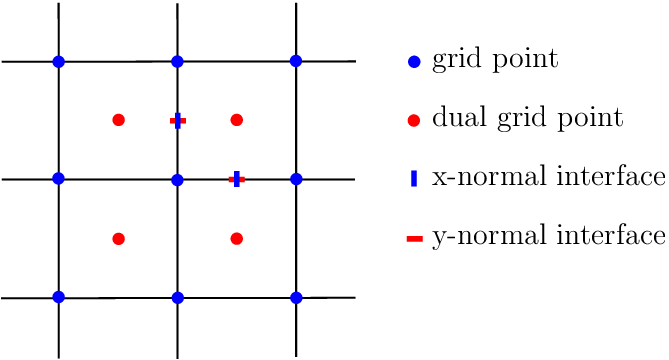}
\caption{\label{1-stencil-2d} Schematic of the dual-grid framework in multi-dimensional FD-CGKS. Blue and red points represent the primary and dual grid nodes, respectively, while the blue and red line segments indicate the corresponding normal interfaces for flux evaluation.}
\end{figure}

In the multidimensional context, space-time coupled methods like the GKS possess an inherent multidimensional nature. This dictates that the numerical flux evaluation depends not only on the flow variables themselves but also on their spatial derivatives along all independent directions \cite{xu2}. Consequently, comprehensive multidimensional flow state information must be provided during the spatial reconstruction at the grid interfaces.
To develop a multidimensional, space-time coupled high-order scheme, this study proposes a novel dual-grid reconstruction method on Cartesian meshes, as illustrated in Fig.~\ref{1-stencil-2d}. In this layout, the blue dots represent the primary grid nodes, while the red dots, located at the centroids of the primary grid virtual cells, are defined as the dual grid points. The dual-grid approach updates the conservative variables on both the primary grid and an identical dual grid, which is offset by half the mesh spacing in each spatial direction.

Based on this dual-grid topology, the conservative variables \( \mathbf{W}_{i,j} \) and \( \mathbf{W}^C_{i,j} \) at the primary and dual grid points, respectively, are updated via Eq.~\eqref{eq-conservation-law-2d}. Meanwhile, the fluxes in the \( x \)- and \( y \)-directions must be evaluated at each interface, as well as at all primary and dual grid points. As illustrated in Fig.~\ref{1-stencil-2d}, the thick blue and red line segments denote the interfaces with different normal directions. The corresponding spatial reconstruction procedures for the flux evaluation are detailed as follows:

\begin{enumerate}
    \item \textbf{Reconstruction at grid interfaces:} Taking the \( x \)-normal interface located on the horizontal primary grid line (denoted by the thick blue segment in Fig.~\ref{1-stencil-2d}) as an example, the reconstruction procedures for other interfaces follow analogously.
    \begin{enumerate}
        \item The left and right states \( \mathbf{W}^{l,r} \) and their normal derivatives \( \mathbf{W}^{l,r}_x \) at this interface are reconstructed using a four-point stencil comprising primary grid points distributed along the \( x \)-direction via a 1D reconstruction (see Section~\ref{sec:reconstruction} and \ref{appendix:C1} and \ref{appendix:C2} for details).
        \item For the tangential derivatives at this interface, a four-point stencil comprising dual grid points distributed along the \( y \)-direction is utilized to perform the 1D reconstruction, yielding \( \mathbf{W}^{C,l,r}_y \). The final reconstructed tangential derivative is obtained by taking the arithmetic average: \( \mathbf{W}^{l}_y = \mathbf{W}^{r}_y = \big( \mathbf{W}^{C,l}_y + \mathbf{W}^{C,r}_y \big)/2 \).
    \end{enumerate}

    \item \textbf{Reconstruction at grid points:} Taking the reconstruction of the \( x \)-derivative at a primary grid point as an illustrative example, the procedures for other derivatives and dual points are carried out similarly.
    \begin{enumerate}
        \item The derivative \( \mathbf{W}_x \) at the primary grid point is reconstructed using a three-point stencil comprising primary grid points distributed along the \( x \)-direction via a 1D reconstruction (see Section~\ref{sec:reconstruction} and \ref{appendix:C3}).
        \item Similarly, the derivative \( \mathbf{W}_y \) at the primary grid point is reconstructed using a three-point stencil comprising primary grid points along the \( y \)-direction, following the identical formulation as in step (a).
    \end{enumerate}
\end{enumerate}

All the aforementioned reconstructions are performed on-the-fly during the flux evaluation, eliminating the need for pre-computation and additional memory storage.
Such a 1D-based reconstruction strategy significantly simplifies the intricate high-order spatial reconstructions typically required by previous space-time coupled GKS on Cartesian grids. Furthermore, the FD discretization inherently avoids the introduction of Gaussian quadrature points along grid interfaces, thereby substantially reducing the computational overhead associated with flux evaluation and spatial reconstruction. These computational benefits constitute the core motivation and methodological foundation of the present study.

\section{Numerical examples}
\label{sec:examples}

This section presents a comprehensive set of benchmark computations to evaluate the performance of the proposed sixth-order FD-CGKS.
First, one- and two-dimensional accuracy tests are conducted to validate the theoretical order of convergence.
Subsequently, a series of one-dimensional problems, including shock tube tests, the interacting blast wave problem with extreme pressure jumps, and the shock-density wave interaction problem, are simulated. These tests are designed to demonstrate the non-oscillatory property of the scheme near strong shocks, its robust shock-capturing capability, as well as its high-resolution and low-dissipation characteristics.
Next, two-dimensional inviscid and viscous flows are investigated, encompassing the double Mach reflection, the Kelvin-Helmholtz instability, the viscous double shear layer, and the viscous shock tube problems. Through these two-dimensional simulations and corresponding quantitative comparisons, the high accuracy and strong robustness of the scheme in resolving complex flow structures are further substantiated.

In the one-dimensional computations, a comparative study between two variants utilizing different nonlinear numerical fluxes, namely FD-CGKS(6-2) and FD-CGKS(6-4), is provided. The numerical results indicate that FD-CGKS(6-4) achieves a superior balance between numerical accuracy and computational robustness. Consequently, for the two-dimensional cases, only the results obtained using the FD-CGKS(6-4) scheme are presented.
The time step is determined by the Courant-Friedrichs-Lewy (CFL) condition, with the CFL number set to \( \text{CFL} \geq 0.4 \) in all test cases, unless otherwise specified.
For viscous flows, the time step is additionally restricted by the viscous stability condition, given by \( \Delta t = \text{CFL} \cdot \frac{\Delta x^2}{8\nu} \), where \( \nu \) is the kinematic viscosity coefficient.

\subsection{Accuracy test}

To verify the spatial accuracy and convergence rate of the proposed FD-CGKS, the advection of a density perturbation is simulated in both one and two dimensions.
For the 1D test, the computational domain is set to \( [0, 2] \) with periodic boundary conditions applied at both ends. The initial condition is given by
\begin{align*}
    \rho(x) = 1 + 0.2\sin(\pi x), \quad U(x) = 1, \quad p(x) = 1, \quad x \in [0,2].
\end{align*}
Similarly, the 2D advection test is conducted on a square domain \( [0, 2] \times [0, 2] \) equipped with periodic boundary conditions in both the \( x \)- and \( y \)-directions. The initial flow field is a natural extension of the 1D case, initialized as
\begin{align*}
    \rho(x,y) = 1 + 0.2\sin\big(\pi (x+y)\big), \quad U(x,y) = 1, \quad V(x,y) = 1, \quad p(x,y) = 1.
\end{align*}

To ensure that the temporal discretization error does not contaminate the assessment of spatial accuracy, the time step is scaled as \( \Delta t \sim \Delta x^2 \) throughout the simulations. The quantitative results of the mesh refinement study are summarized in Table~\ref{1d-accuracy} for the 1D case and Table~\ref{2d-accuracy} for the 2D case.
Only the results of FD-CGKS(6-2) are presented in these cases for illustrative purposes, while FD-CGKS(6-4) is also capable of reaching the theoretical sixth-order convergence.
These tables present the \( L^1 \) and \( L^\infty \) errors of the density field, along with their corresponding convergence orders. As the mesh is progressively refined, the proposed scheme consistently achieves the designed sixth-order spatial accuracy in both the one- and two-dimensional simulations.

\begin{table}
	\begin{center}
		\def\temptablewidth{0.90\textwidth}
		{\rule{\temptablewidth}{0.70pt}}
		\begin{tabular*}{\temptablewidth}{@{\extracolsep{\fill}}c|cc|cc}
			
			mesh size & $L^1$ error & Order & $L^{\infty}$ error & Order  \\
			\hline
            1/4  & 4.8840e-05 &       & 7.9701e-05 &      \\
            1/8  & 8.9234e-07 & 5.77  & 1.3542e-06 & 5.88 \\
            1/16 & 1.4734e-08 & 5.92  & 2.2934e-08 & 5.88 \\
            1/32 & 2.3760e-10 & 5.95  & 3.7096e-10 & 5.95 \\
            1/64 & 3.7639e-12 & 5.98  & 5.9014e-12 & 5.97 \\
		\end{tabular*}
		{\rule{\temptablewidth}{0.7pt}}
	\end{center}
	\vspace{-1mm} \caption{\label{1d-accuracy} 1D accuracy test: the \( L^1 \) and \( L^\infty \) errors and convergence orders of the 6th-order compact FD-CGKS(6-2).}
\end{table}

\begin{table}
	\begin{center}
		\def\temptablewidth{0.90\textwidth}
		{\rule{\temptablewidth}{0.70pt}}
		\begin{tabular*}{\temptablewidth}{@{\extracolsep{\fill}}c|cc|cc}
			
			mesh size & $L^1$ error & Order & $L^{\infty}$ error & Order  \\
			\hline
            2/5  & 3.7278e-03 &       & 1.4016e-03 &      \\
            1/5  & 5.1897e-05 & 6.17  & 2.0373e-05 & 6.10 \\
            1/10 & 8.8195e-07 & 5.88  & 3.4686e-07 & 5.88 \\
            1/20 & 1.4870e-08 & 5.89  & 5.8351e-09 & 5.89 \\
            1/40 & 2.4184e-10 & 5.94  & 9.5053e-11 & 5.94 \\
		\end{tabular*}
		{\rule{\temptablewidth}{0.7pt}}
	\end{center}
	\vspace{-1mm} \caption{\label{2d-accuracy} 2D accuracy test: the \( L^1 \) and \( L^\infty \) errors and convergence orders of the 6th-order compact FD-CGKS(6-2).}
\end{table}

\subsection{One-dimensional shock tube problems}

To evaluate the essentially non-oscillatory property of the proposed schemes in the presence of shock waves, as well as their resolution capabilities for contact discontinuities, two classic 1D Riemann problems, namely the Sod and Lax shock tube problems, are simulated.
For the Sod problem \cite{sod1978_test}, the initial condition is given by
\begin{align*}
(\rho, U, p) =
\begin{cases}
(1, 0, 1), & 0 \leq x < 0.5, \\
(0.125, 0, 0.1), & 0.5 \leq x \leq 1.
\end{cases}
\end{align*}
For the Lax problem \cite{lax1954_test}, the initial condition is specified as
\begin{align*}
(\rho, U, p) =
\begin{cases}
(0.445, 0.698, 3.528), & 0 \leq x < 0.5, \\
(0.5, 0, 0.571), & 0.5 \leq x \leq 1.
\end{cases}
\end{align*}
For both test cases, the computational domain is set to \( [0, 1] \) and is discretized using a uniform mesh with \( N = 100 \) grid points.
Zero-gradient boundary conditions are applied at both ends. The output times are $t=0.2$ for the Sod problem and $t=0.16$ for the Lax problem.

\begin{figure}[!htb]
\centering
\includegraphics[width=0.495\textwidth]{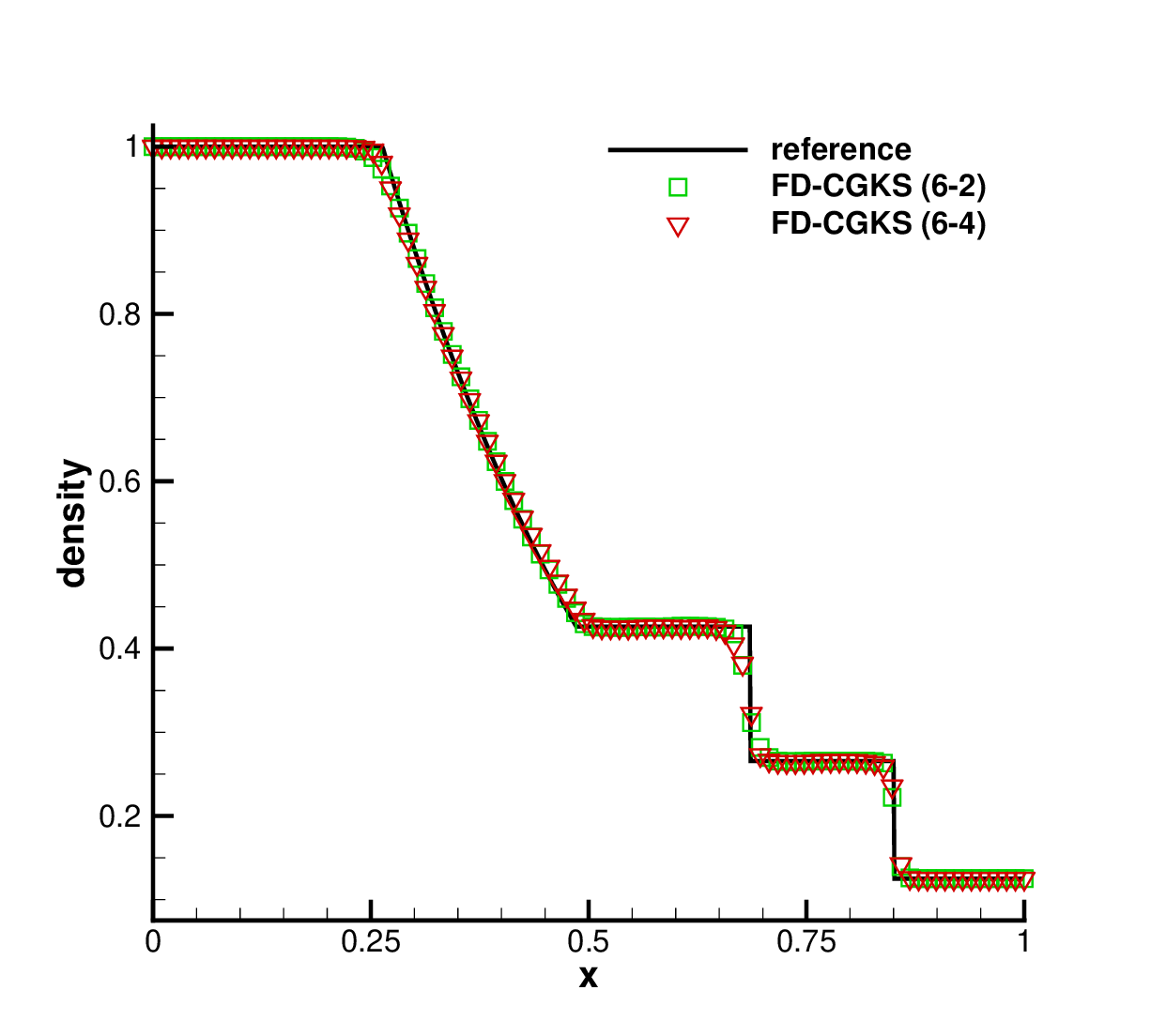}
\includegraphics[width=0.495\textwidth]{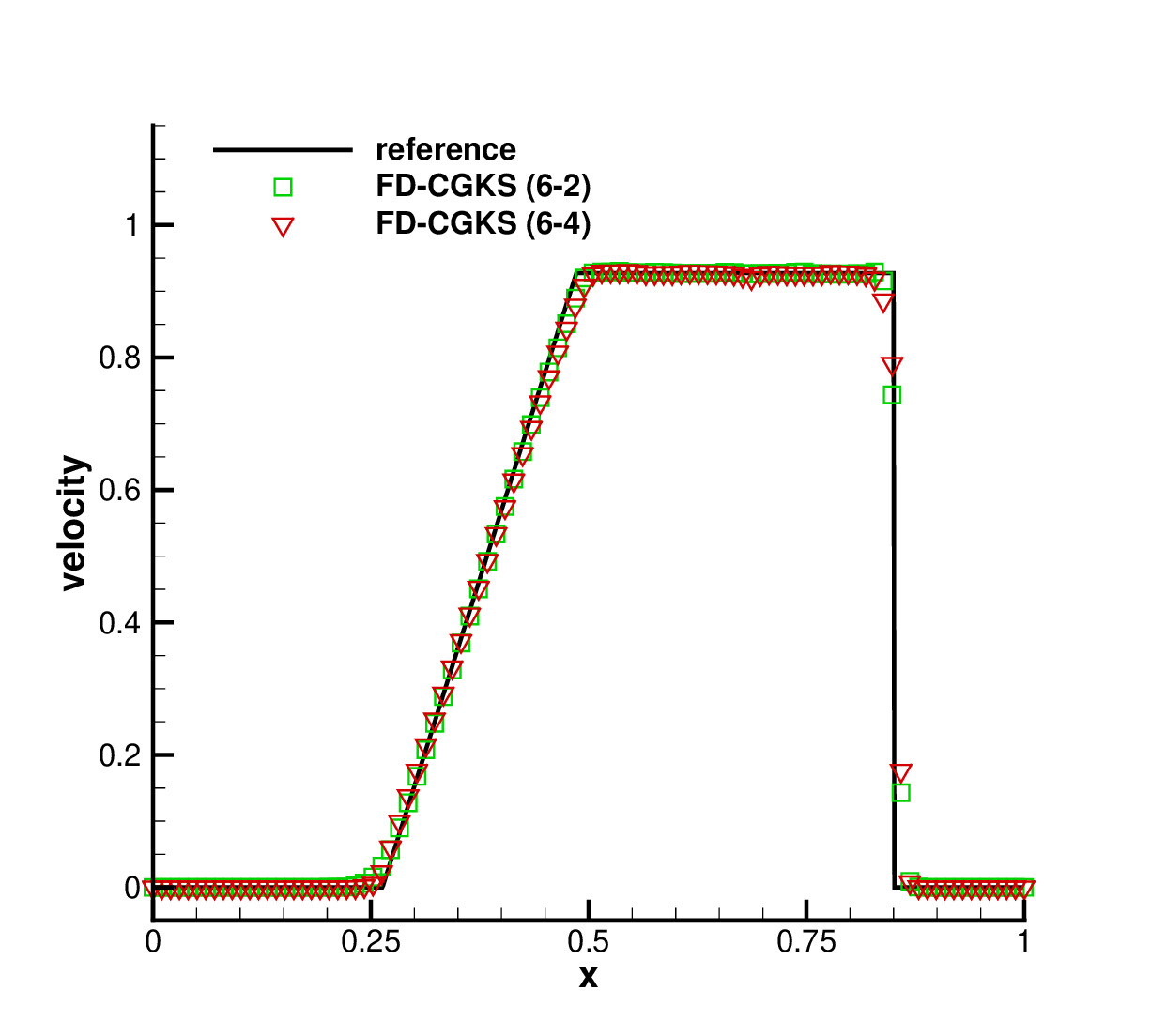}
\caption{\label{1d-sod} 1D Sod shock tube problem: The density and velocity distributions are computed using the FD-CGKS(6-4) and FD-CGKS(6-2) schemes with \( 100 \) uniform mesh points. }
\end{figure}

\begin{figure}[!htb]
\centering
\includegraphics[width=0.495\textwidth]{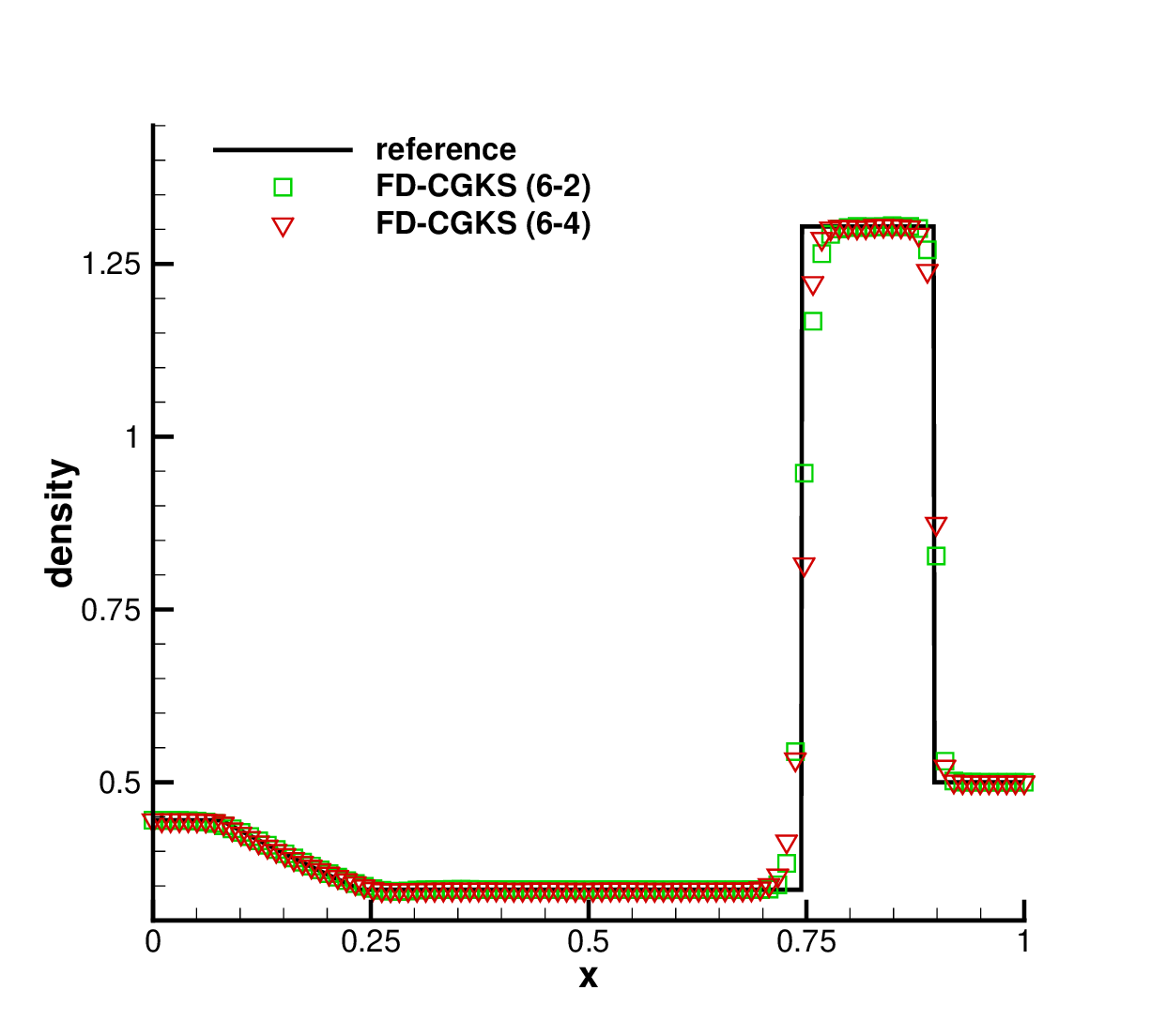}
\includegraphics[width=0.495\textwidth]{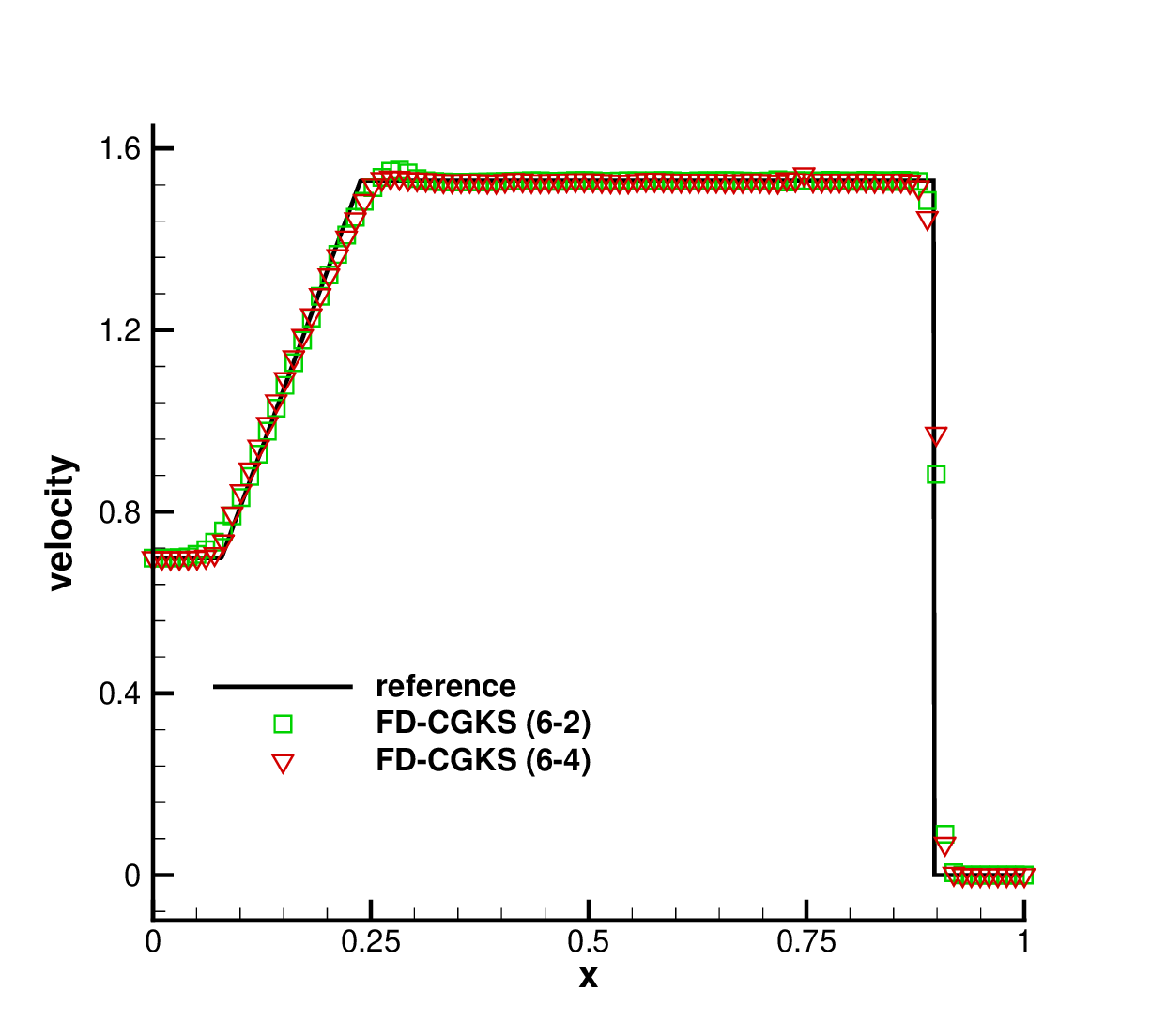}\\
\caption{\label{1d-lax} 1D Lax shock tube problem: The density and velocity distributions are computed using the FD-CGKS(6-4) and FD-CGKS(6-2) schemes with \( 100 \) uniform mesh points. Both schemes demonstrate high resolution and essentially non-oscillatory behavior.}
\end{figure}

The numerical simulations are conducted using the FD-CGKS equipped with two different nonlinear numerical fluxes, denoted as FD-CGKS(6-4) and FD-CGKS(6-2). The computed density and velocity distributions for the Sod and Lax problems are presented in Fig.~\ref{1d-sod} and Fig.~\ref{1d-lax}, respectively. The exact solutions are plotted as reference lines.
As observed from the results, both schemes successfully achieve a sharp, essentially non-oscillatory capturing of the shock waves and contact discontinuities while maintaining high resolution across these discontinuities. Furthermore, the comparisons demonstrate that the FD-CGKS(6-4) and FD-CGKS(6-2) schemes exhibit nearly identical resolving capabilities for these benchmark problems.

\subsection{One-dimensional blast wave problems}

Two test cases featuring extreme initial pressure jumps are employed to validate the robustness of the newly proposed FD-CGKS in handling such severe conditions and capturing strong shock waves.
The first test case is the Woodward-Colella blast wave problem \cite{Case-Woodward}. The initial condition is given by
\begin{align*}
(\rho, U, p) =
\begin{cases}
(1, 0, 1000), & 0 \leq x < 0.1, \\
(1, 0, 0.01), & 0.1 \leq x < 0.9, \\
(1, 0, 100),  & 0.9 \leq x \leq 1.
\end{cases}
\end{align*}
The computational domain is \( [0, 1] \), and reflecting boundary conditions are imposed at both ends of the domain. The domain is discretized using a uniform mesh with \( N = 400 \) nodes. The simulation is advanced to a final computation time of \( t = 0.038 \).

The second test case is the 1D Sedov blast wave problem \cite{xiong2016-test}. The initial condition is characterized by a uniform density of \( \rho = 1 \) and a zero velocity field \( U = 0 \) throughout the domain. The initial total energy is set to \( 10^{-8} \) everywhere, except in the center grid where it is specified as a constant \( E_0 / \Delta x \), with \( E_0 = 3.2 \times 10^6 \). The computational domain is discretized using a uniform mesh with \( N = 800 \) nodes. The simulation is run up to a final time of \( t = 0.001 \).

\begin{figure}[!htb]
\centering
\includegraphics[width=0.495\textwidth]{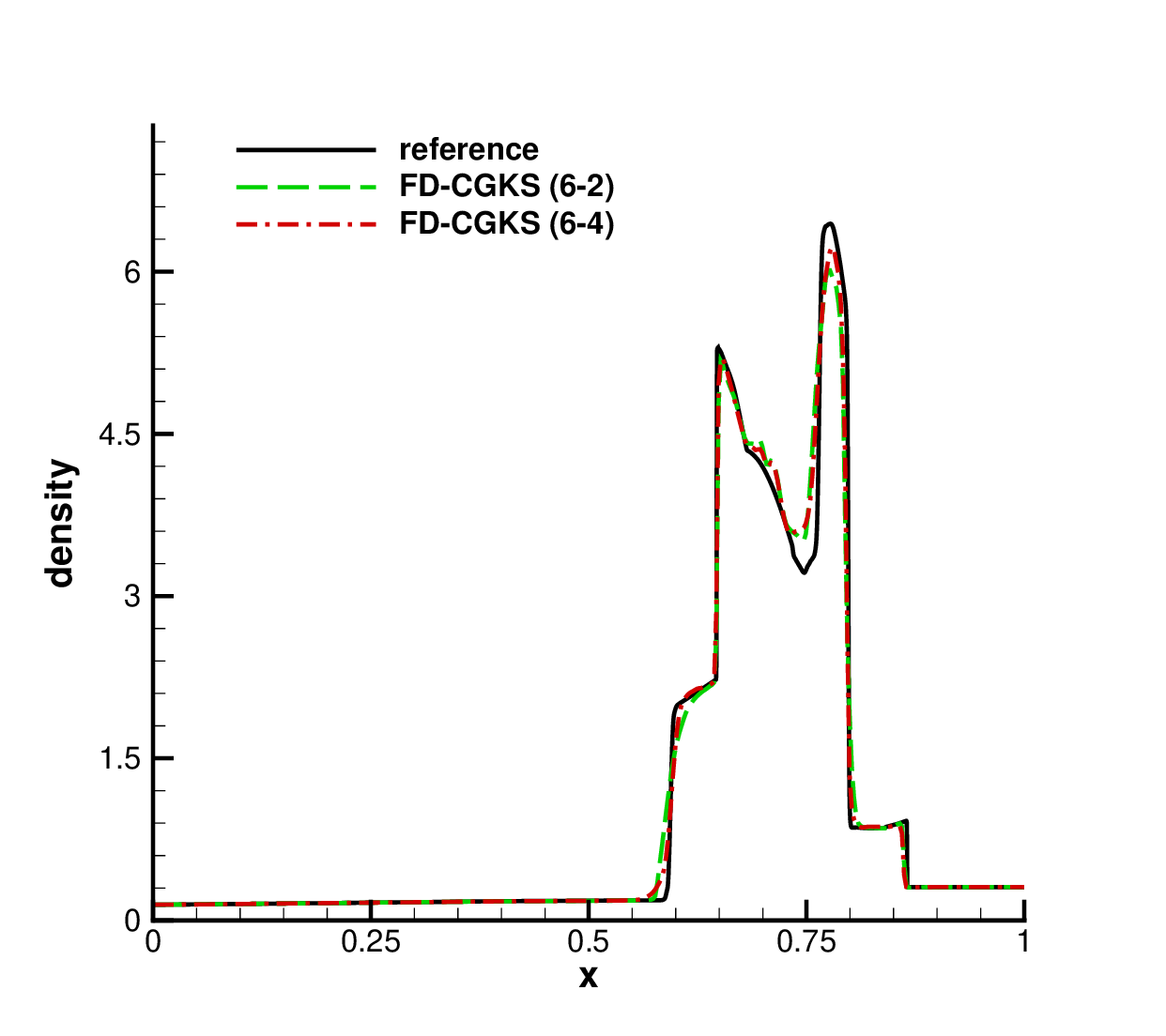}
\includegraphics[width=0.495\textwidth]{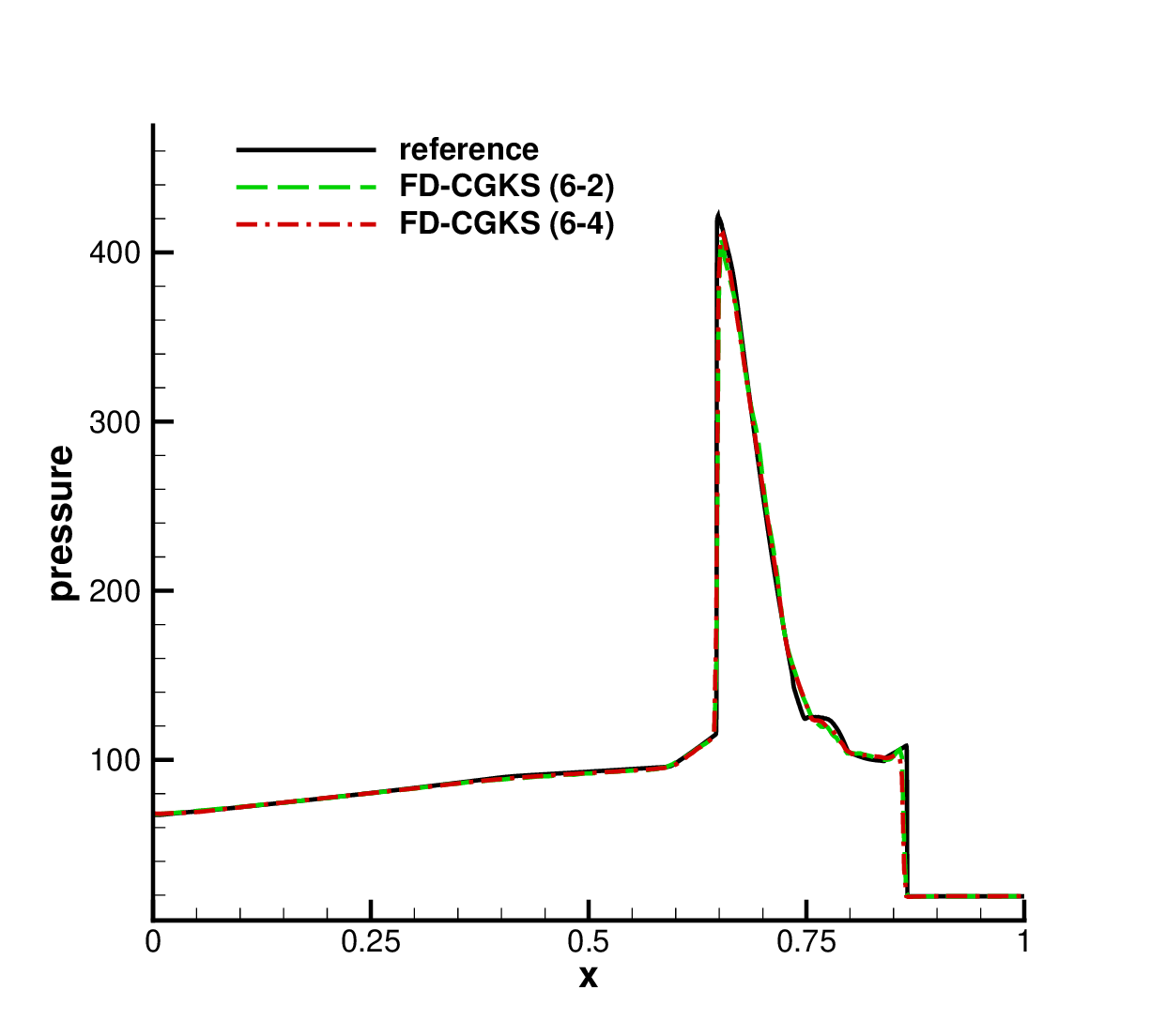}\\
\caption{\label{1d-blast} Woodward-Colella blast wave problem: The density and pressure distributions are computed using the FD-CGKS(6-4) and FD-CGKS(6-2) with \( 400 \) uniform mesh points. Both nonlinear numerical fluxes exhibit robust capturing of the strong interacting shock waves. }
\end{figure}

\begin{figure}[!htb]
\centering
\includegraphics[width=0.495\textwidth]{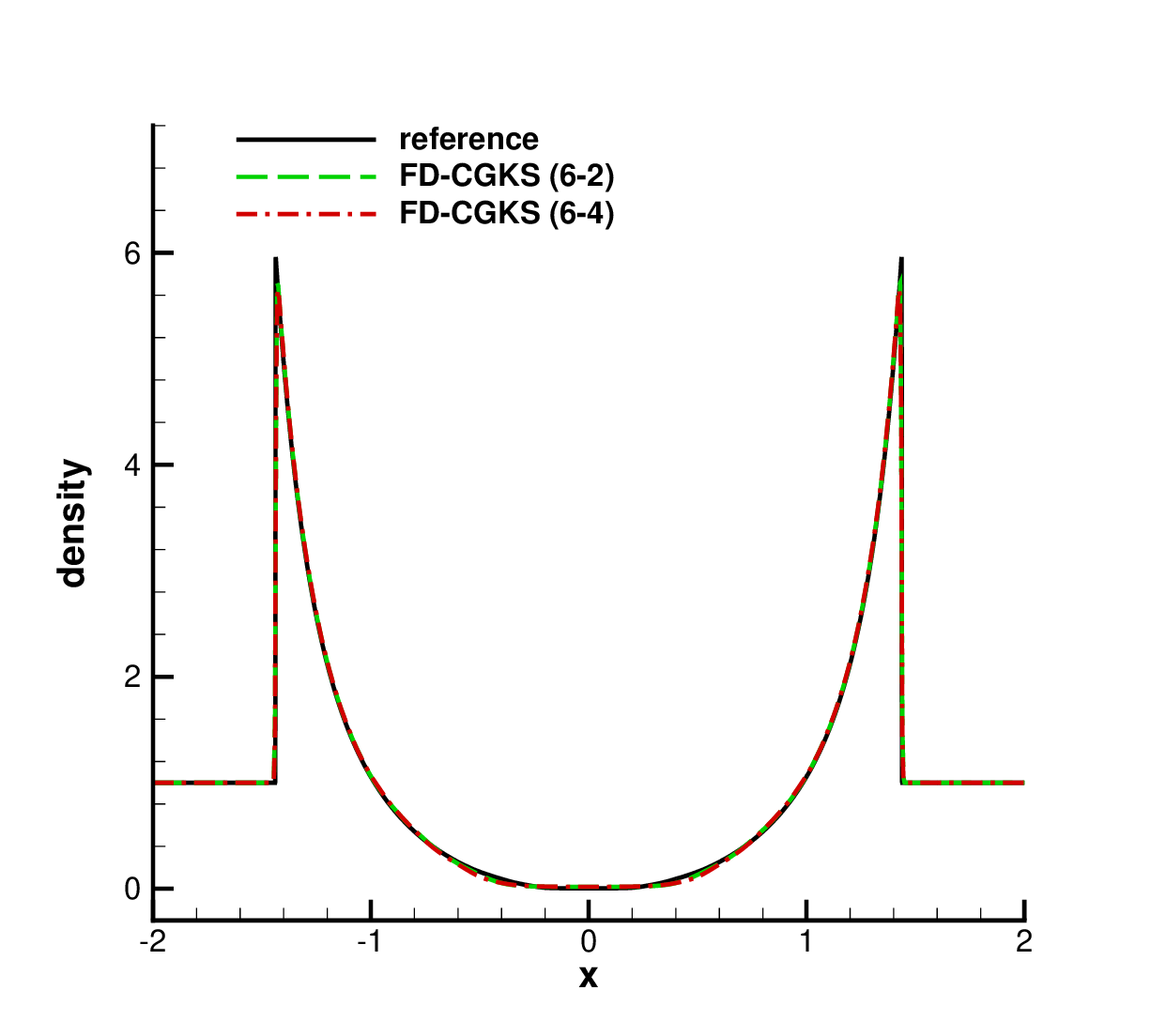}
\includegraphics[width=0.495\textwidth]{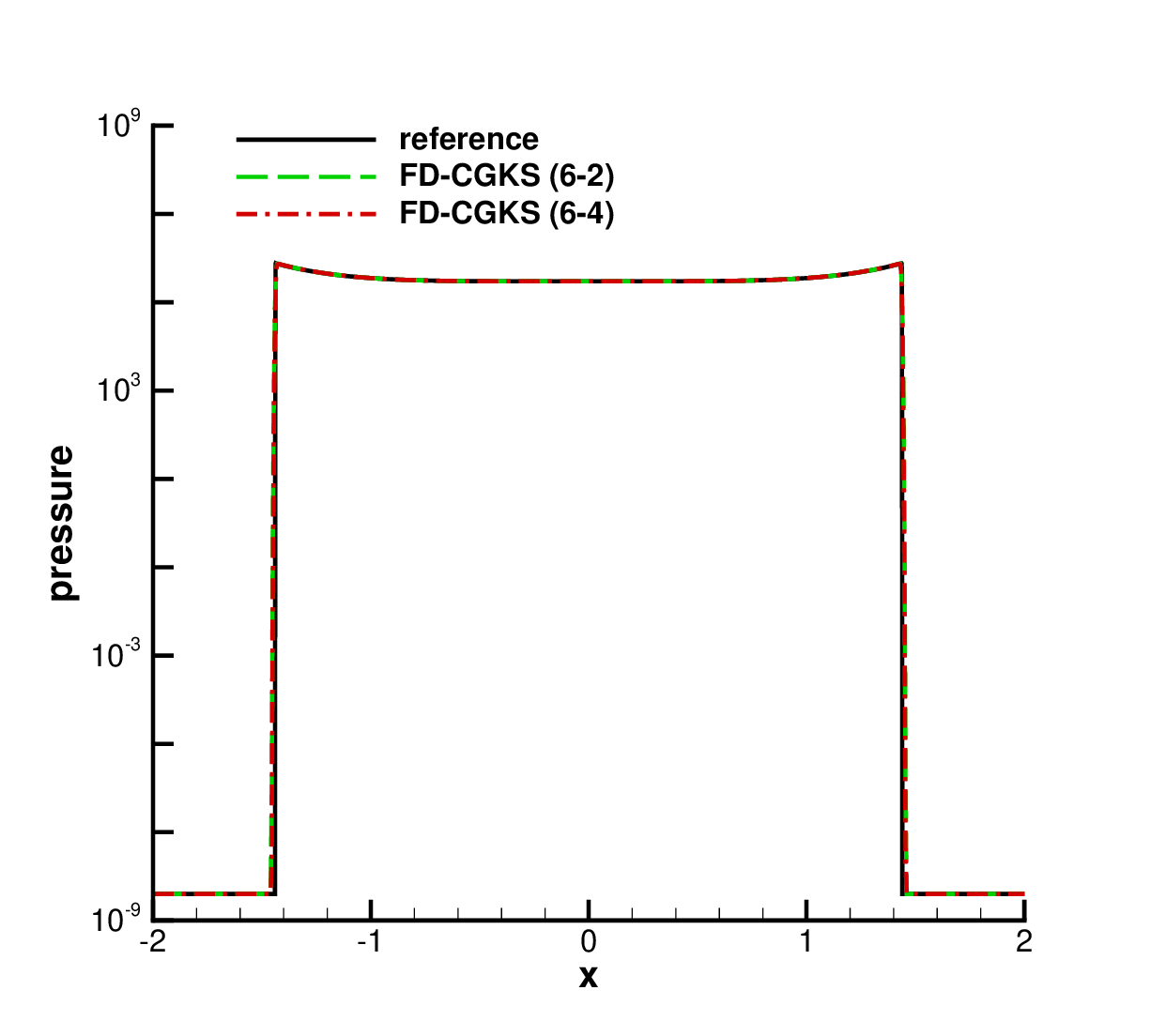}\\
\caption{\label{1d-sedov} Sedov blast wave problem: The density and pressure distributions are computed using the FD-CGKS(6-4) and FD-CGKS(6-2) with \( 800 \) uniform mesh points. The results demonstrate the capability of the schemes to handle extreme initial energy concentrations and strong pressure discontinuities. }
\end{figure}

Fig.~\ref{1d-blast} and Fig.~\ref{1d-sedov} present the computed density and pressure distributions for the two blast wave problems obtained via the proposed sixth-order FD-CGKS. These results demonstrate the scheme's remarkable robustness in resolving extremely strong discontinuities and shock waves without spurious oscillations.
While the two variants equipped with different nonlinear numerical fluxes, FD-CGKS(6-2) and FD-CGKS(6-4), yield highly consistent results, FD-CGKS(6-4) provides sharper resolution for contact discontinuities in the Woodward-Colella test. Additionally, because these extreme cases involve highly unsteady and strong shock waves, the application of the nonlinear S2O4 time discretization method is indispensable for preserving numerical stability.

\subsection{1D Titarev-Toro problems}

To assess the high-resolution capabilities of the proposed FD-CGKS, the Titarev-Toro problem \cite{titarev2004finite} is considered. Characterized by the interaction between a shock wave and high-wavenumber density perturbations, this test case serves as a stringent benchmark to demonstrate the scheme's resolving power and low numerical dissipation.
The initial condition for this problem is given by
\begin{align*}
(\rho, U, p) =
\begin{cases}
(1.515695, 0.523346, 1.805), & x \leq -4, \\
(1 + 0.1\sin(20\pi x), 0, 1), & x > -4.
\end{cases}
\end{align*}
The computational domain is defined as \( [-5, 5] \). An inflow boundary condition is imposed at the left boundary, while a fixed wave profile is specified at the right boundary.
The simulation is performed on a uniform mesh with \( N = 1000 \) grid points and is advanced to a final computation time of \( t = 5.0 \).

\begin{figure}[!htb]
\centering
\includegraphics[width=0.495\textwidth]{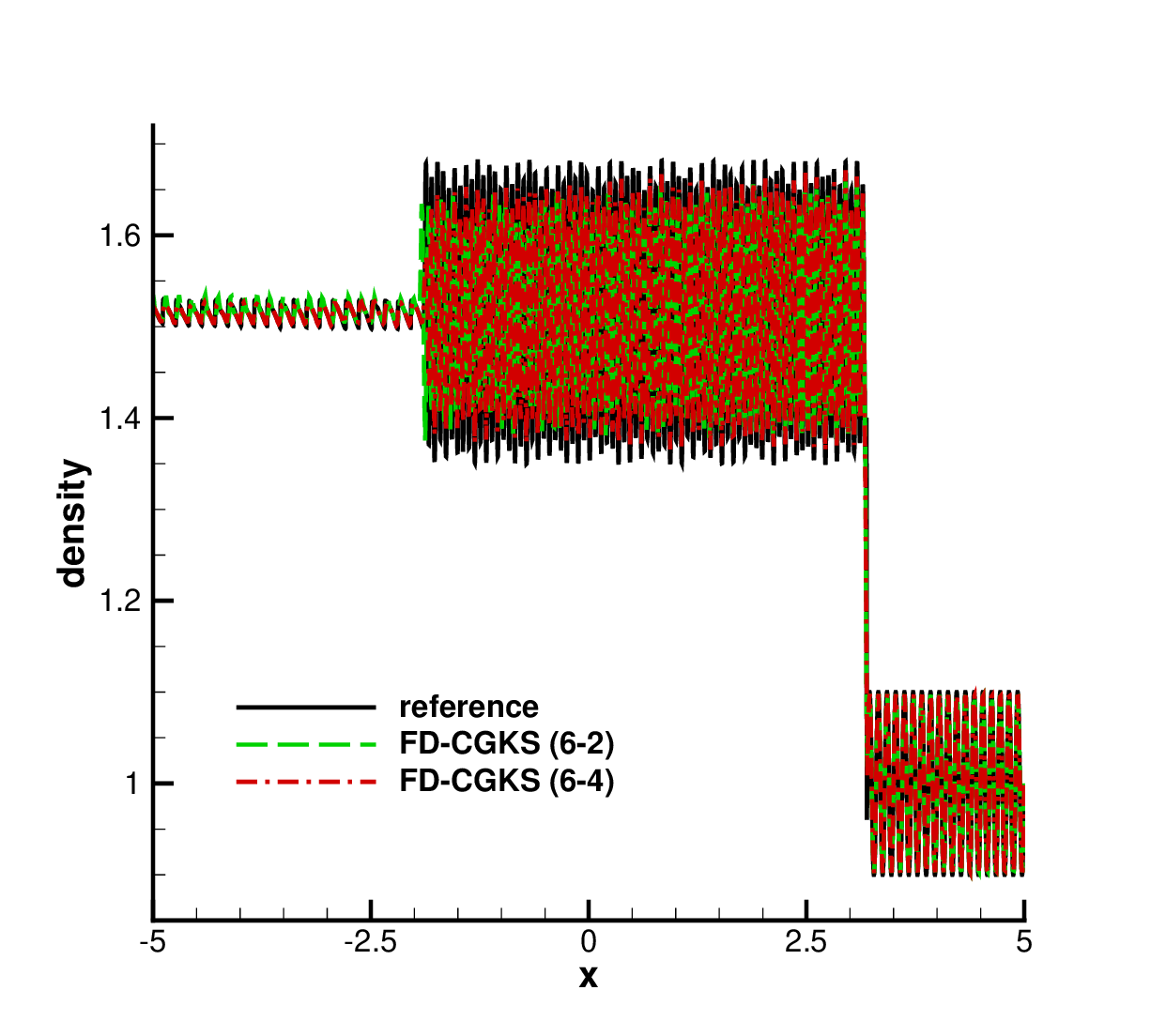}
\includegraphics[width=0.495\textwidth]{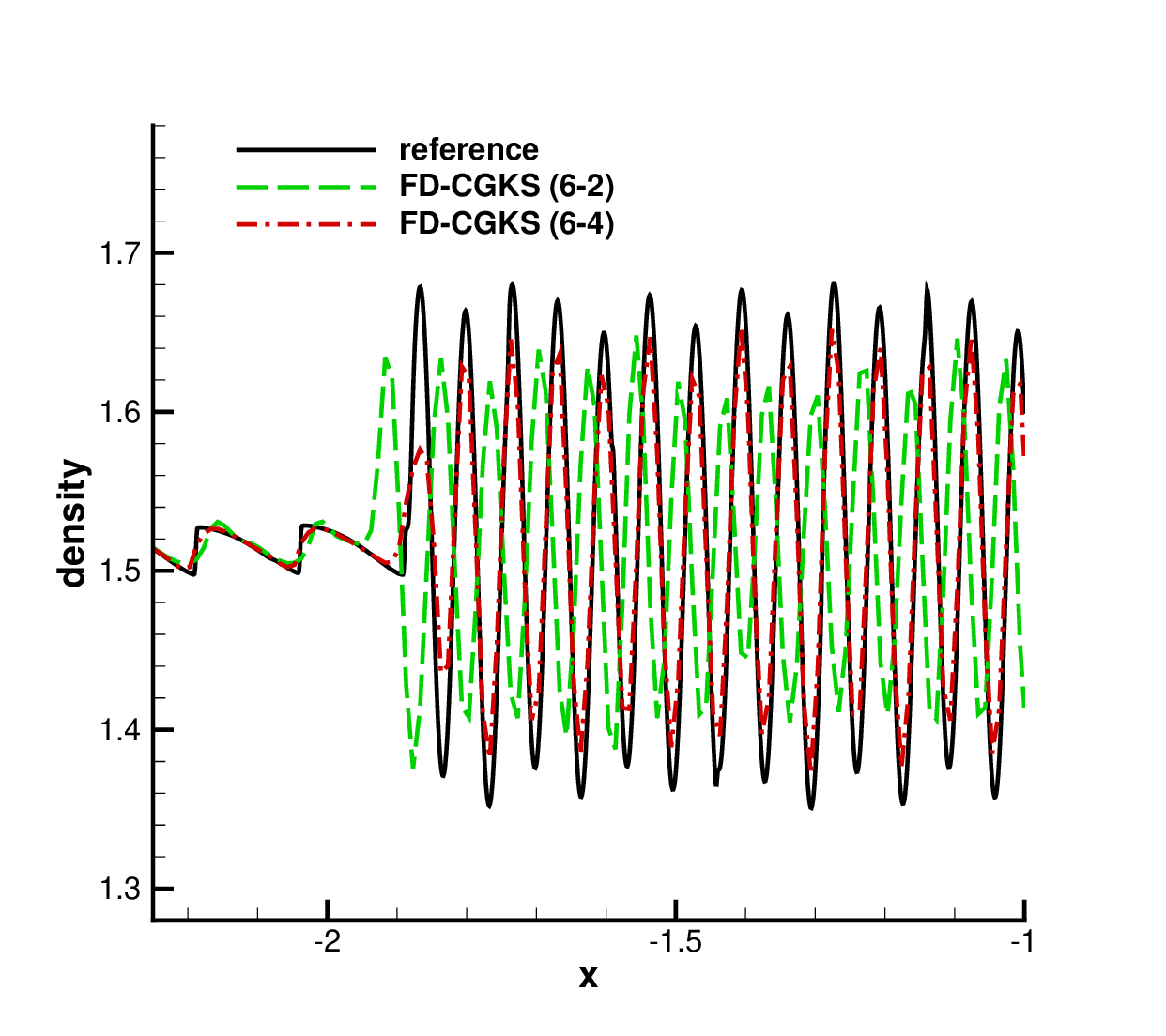}\\
\caption{\label{1d-toro} Titarev-Toro problem: The density distribution is computed using the sixth-order FD-CGKS with two different nonlinear fluxes, FD-CGKS(6-4) and FD-CGKS(6-2), on a uniform mesh of \( 1000 \) points. The right panel displays an enlarged view of the high-wavenumber region, demonstrating the superior phase accuracy and resolving power of the FD-CGKS(6-4) scheme.}
\end{figure}

Fig.~\ref{1d-toro} presents the numerical results obtained using the sixth-order FD-CGKS equipped with the two different nonlinear numerical fluxes. The right panel of Fig.~\ref{1d-toro} provides a locally magnified view of the oscillatory region to better illustrate the flow details.
As observed from the comparisons, the FD-CGKS(6-4) scheme, which employs a higher-order nonlinear flux, achieves significantly higher accuracy. In particular, it captures the phase and amplitude of the high-wavenumber density waves much more accurately than the FD-CGKS(6-2) variant, demonstrating the excellent balance between strong robustness and high accuracy of utilizing higher-order spatial difference discretization in regions with complex wave structures.

\begin{figure}[!htb]
\centering
\includegraphics[width=0.80\textwidth]{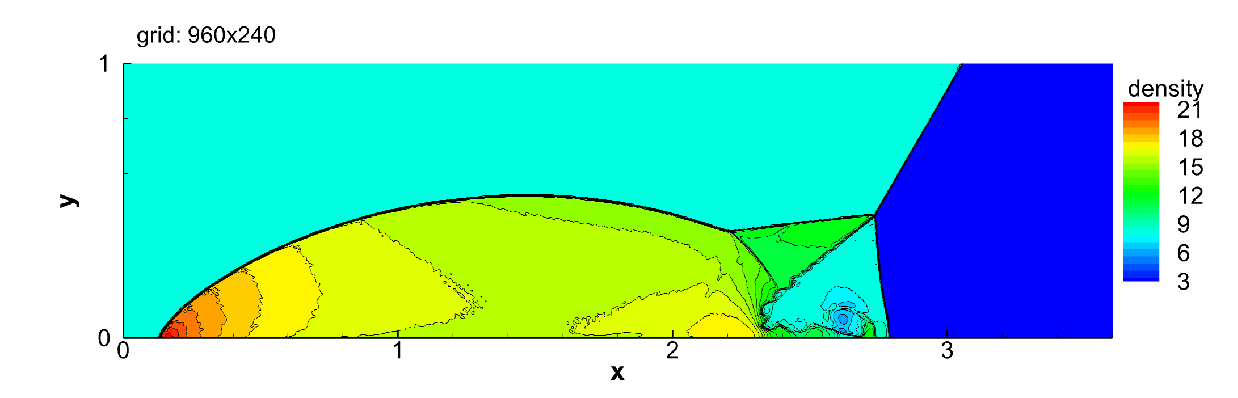}
\includegraphics[width=0.80\textwidth]{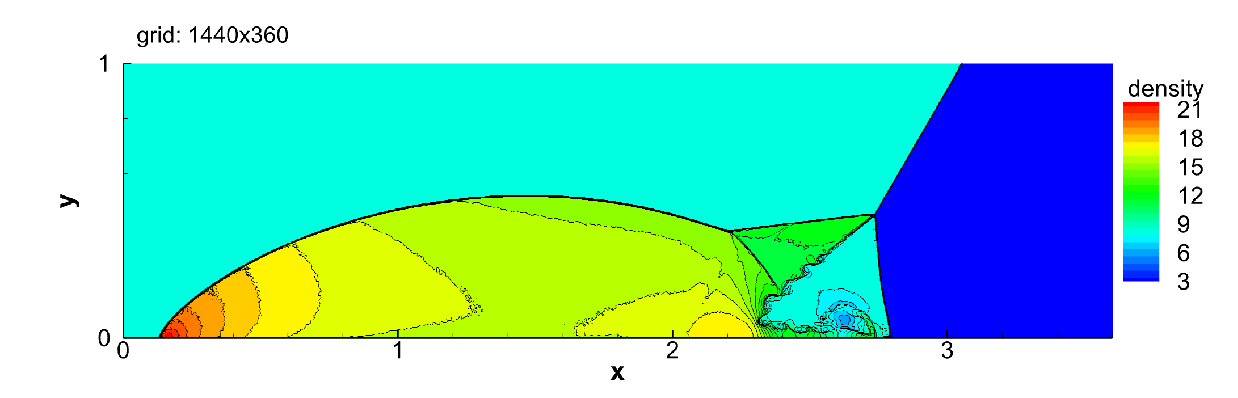}
\caption{\label{2d-doublemach-1} Double Mach reflection problem: Density contours on uniform grids of \( 960\times 240 \) (top) and \( 1440\times 360 \) (bottom) points. The scheme robustly captures the strong shock waves without noticeable spurious oscillations.}
\end{figure}

\begin{figure}[!htb]
\centering
\includegraphics[width=0.45\textwidth]{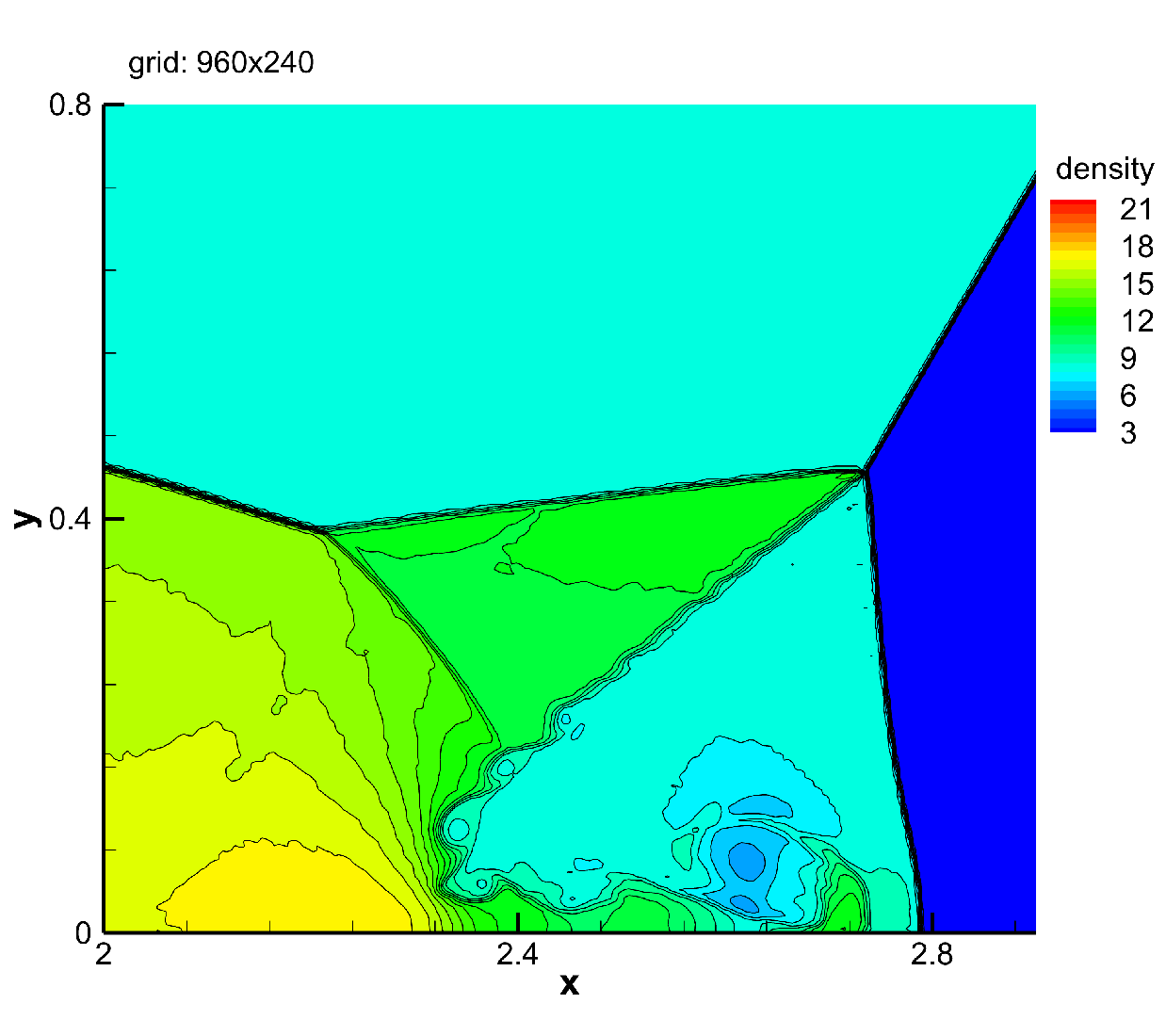}
\includegraphics[width=0.45\textwidth]{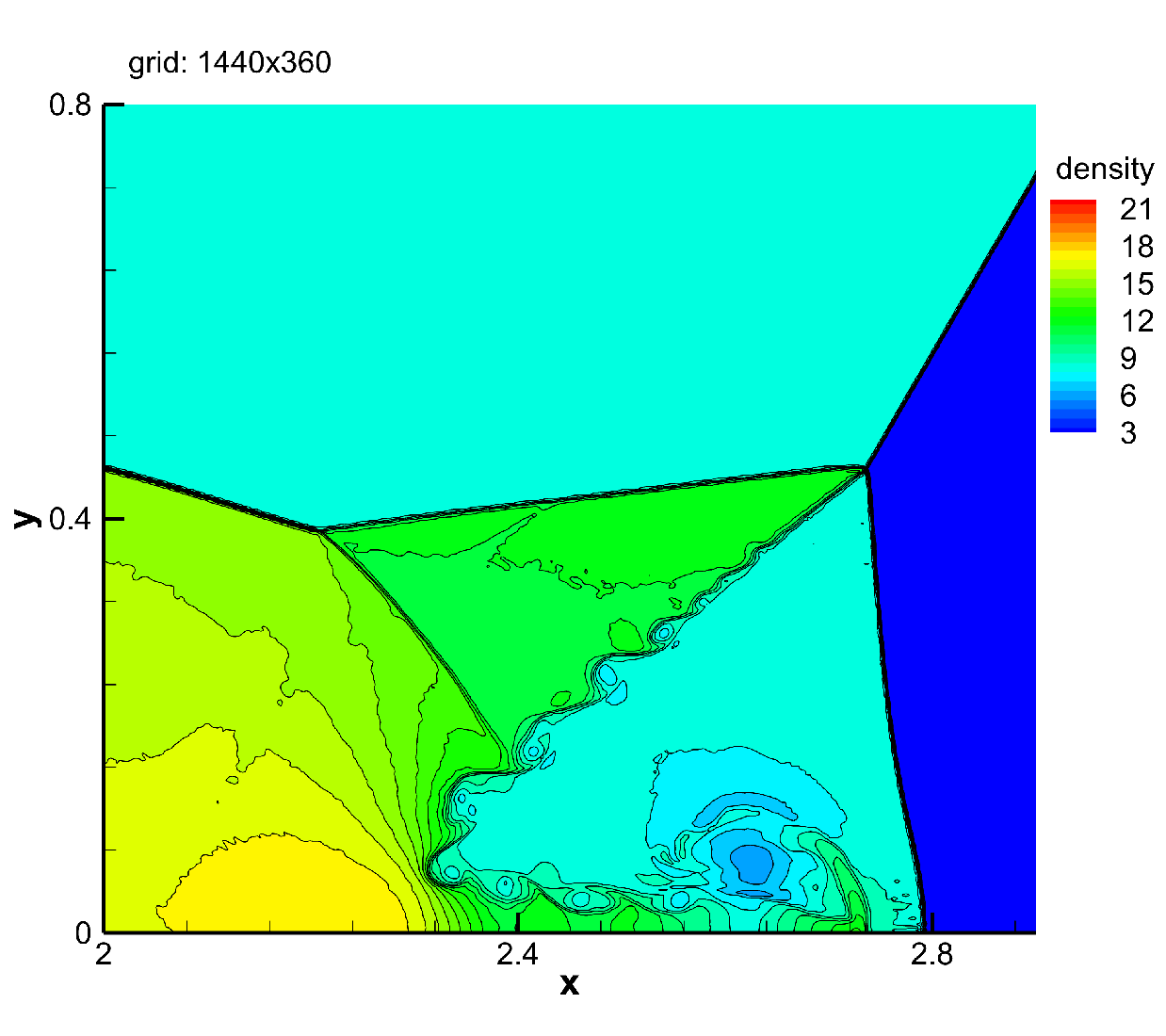}
\caption{\label{2d-doublemach-2} Double Mach reflection problem: Close-up views of the density contours corresponding to the grids of \( 960\times 240 \) (left) and \( 1440\times 360 \) (right) points. The scheme sharply resolves the Mach stems and the Kelvin-Helmholtz vortex roll-ups along the slip lines.}
\end{figure}

\subsection{Double Mach reflection problem}

The two-dimensional inviscid double Mach reflection problem \cite{Case-Woodward} is simulated to evaluate the performance of the proposed scheme in resolving complex shock-shock interactions and small-scale vortical structures.
The computational domain is defined as \( [0, 4] \times [0, 1] \). A reflecting solid wall is located along the bottom boundary, starting from \( x = 1/6 \) and extending to \( x = 4 \). Initially, a right-moving Mach \( 10 \) shock is positioned at \( (x, y) = (1/6, 0) \), inclined at an angle of \( 60^\circ \) relative to the \( x \)-axis. The initial states for the post-shock and pre-shock regions are given respectively by
\begin{align*}
(\rho, U, V, p) =
\begin{cases}
(8, 4.125\sqrt{3}, -4.125,116.5), ~~& \mathrm{post-shock}, \\
(1.4, 0, 0, 1), ~~& \mathrm{pre-shock}.
\end{cases}
\end{align*}
Reflecting boundary conditions are applied along the solid wall, while the outflow conditions are imposed on the remainder of the bottom boundary (\( 0 \leq x < 1/6 \)). At the top boundary, the flow variables are dynamically set to track the exact motion of the Mach \( 10 \) shock. Inflow and outflow boundary conditions are applied at the left and right boundaries, respectively.

The simulations are performed on uniform meshes of \( 960\times 240 \) and \( 1440\times 360 \) grid points. The density contours at \( t = 0.2 \) are presented in Fig.~\ref{2d-doublemach-1}, while the corresponding close-up views of the Mach stem region are shown in Fig.~\ref{2d-doublemach-2}.
The present FD-CGKS(6-4) suppresses spurious oscillations near the strong shock discontinuities, demonstrating its robust shock-capturing capability.
Moreover, the intricate roll-up of the slip lines emanating from the primary and secondary triple points, driven by the Kelvin-Helmholtz instability, is clearly captured. The sharp resolution of these small-scale vortices confirms the high resolution and low numerical dissipation of the proposed scheme.

\subsection{Kelvin--Helmholtz instability problem}

The Kelvin--Helmholtz (KH) instability, driven by the relative motion between two fluid layers or a strong velocity gradient within a single fluid, is a classical benchmark for evaluating the resolving capability of numerical schemes in capturing the nonlinear evolution of shear-driven vortical structures.
In this test, we adopt the stratified KH instability problem \cite{san2015-KHI} to assess the ability of the present sixth-order FD-CGKS in resolving multi-scale flow features arising from the nonlinear growth of interfacial perturbations.
The computational domain is $\Omega = [0,1]\times[0,1]$, and the initial flow field is divided into three zones according to
\begin{align*}
  (\rho,\, U,\, V,\, p) =
  \begin{cases}
    (1,\; 0.5,\; V_0,\; 2.5), & 0 < y < 0.25, \\[4pt]
    (2,\;-0.5,\; V_0,\; 2.5), & 0.25 \le y \le 0.75, \\[4pt]
    (1,\; 0.5,\; V_0,\; 2.5), & 0.75 < y < 1,
  \end{cases}
\end{align*}
where $V_0 = 0.01\sin(2\pi x/L)$ with $L=1$ is a first-mode sinusoidal perturbation imposed on the transverse velocity to trigger the instability.
The perturbation initiates a nonlinear roll-up of the shear layer, which subsequently develops into increasingly complex vortical structures and eventually transitions toward a fully turbulent state. Periodic boundary conditions are applied in both directions.

\begin{figure}[!htb]
\centering
\includegraphics[width=0.325\textwidth]{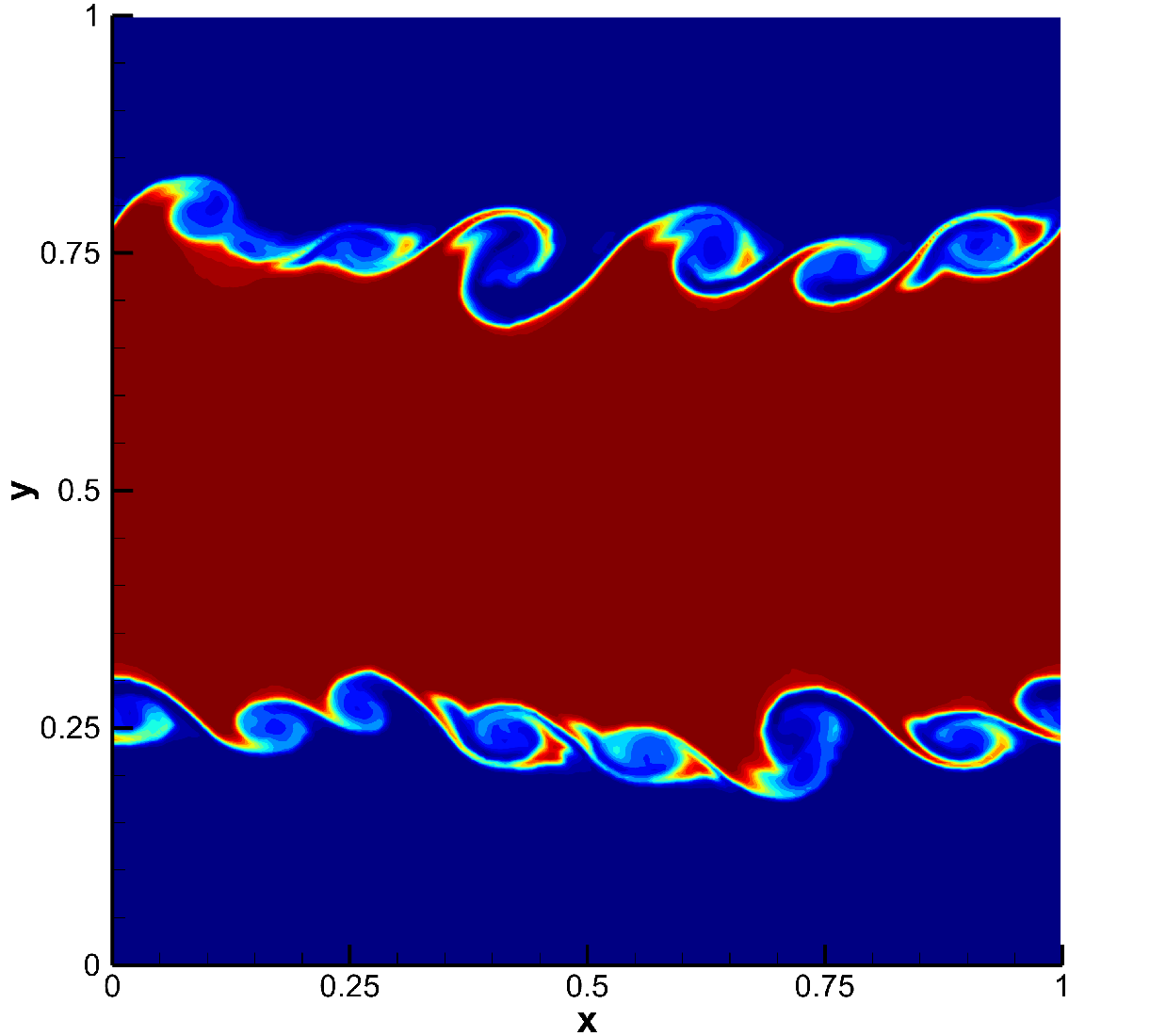}
\includegraphics[width=0.325\textwidth]{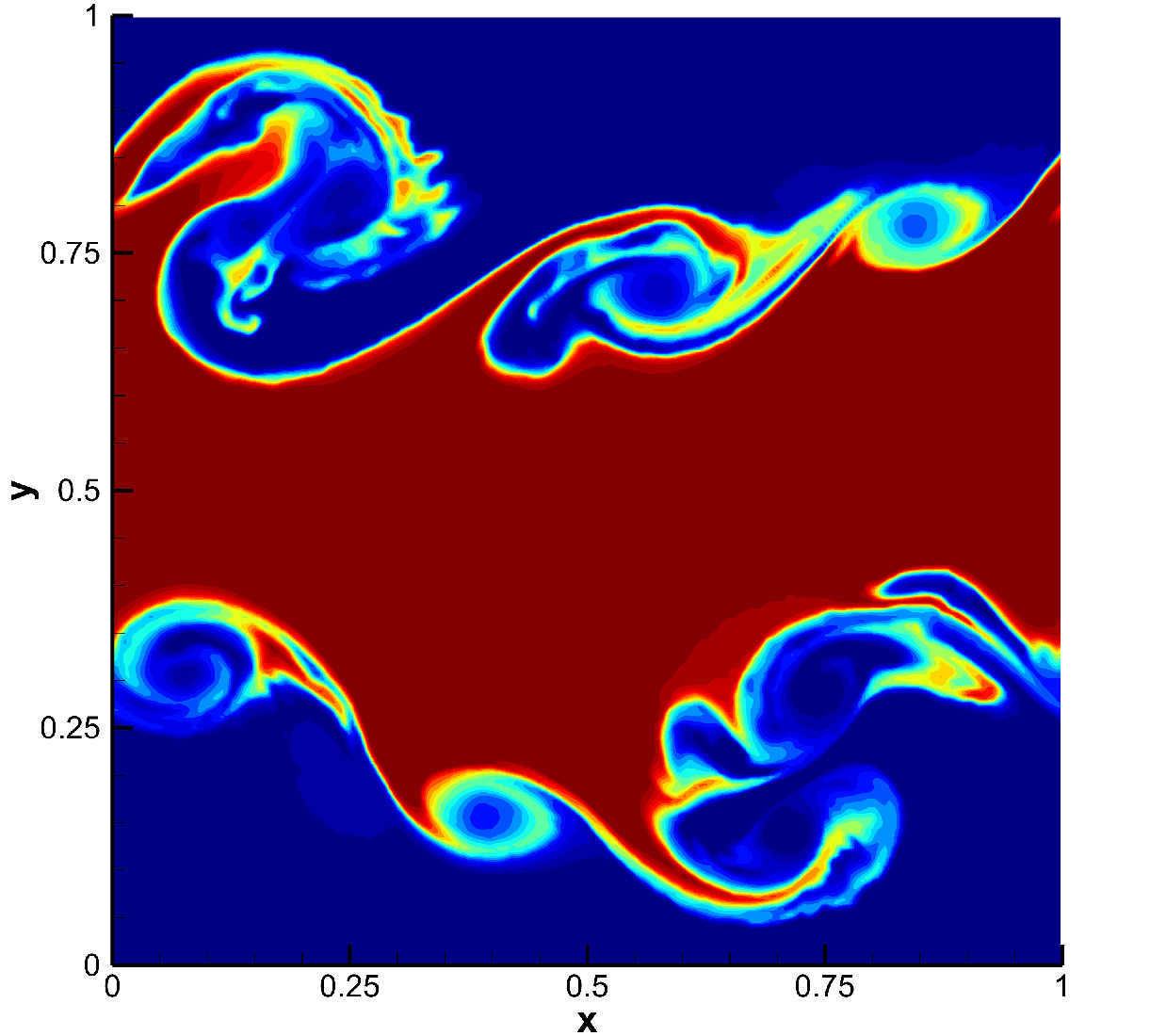}
\includegraphics[width=0.325\textwidth]{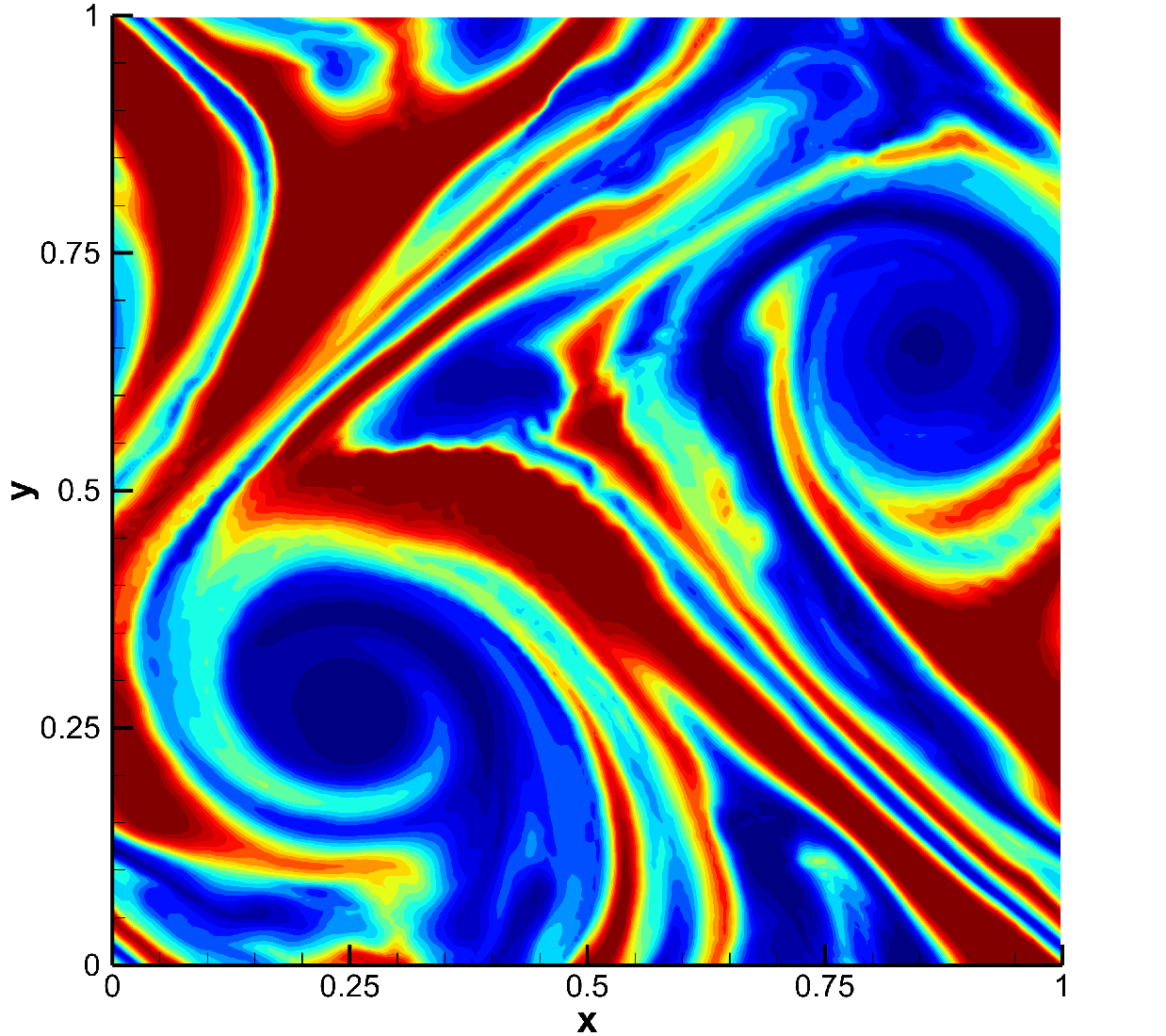}\\
\vspace{0.25cm}
\includegraphics[width=0.325\textwidth]{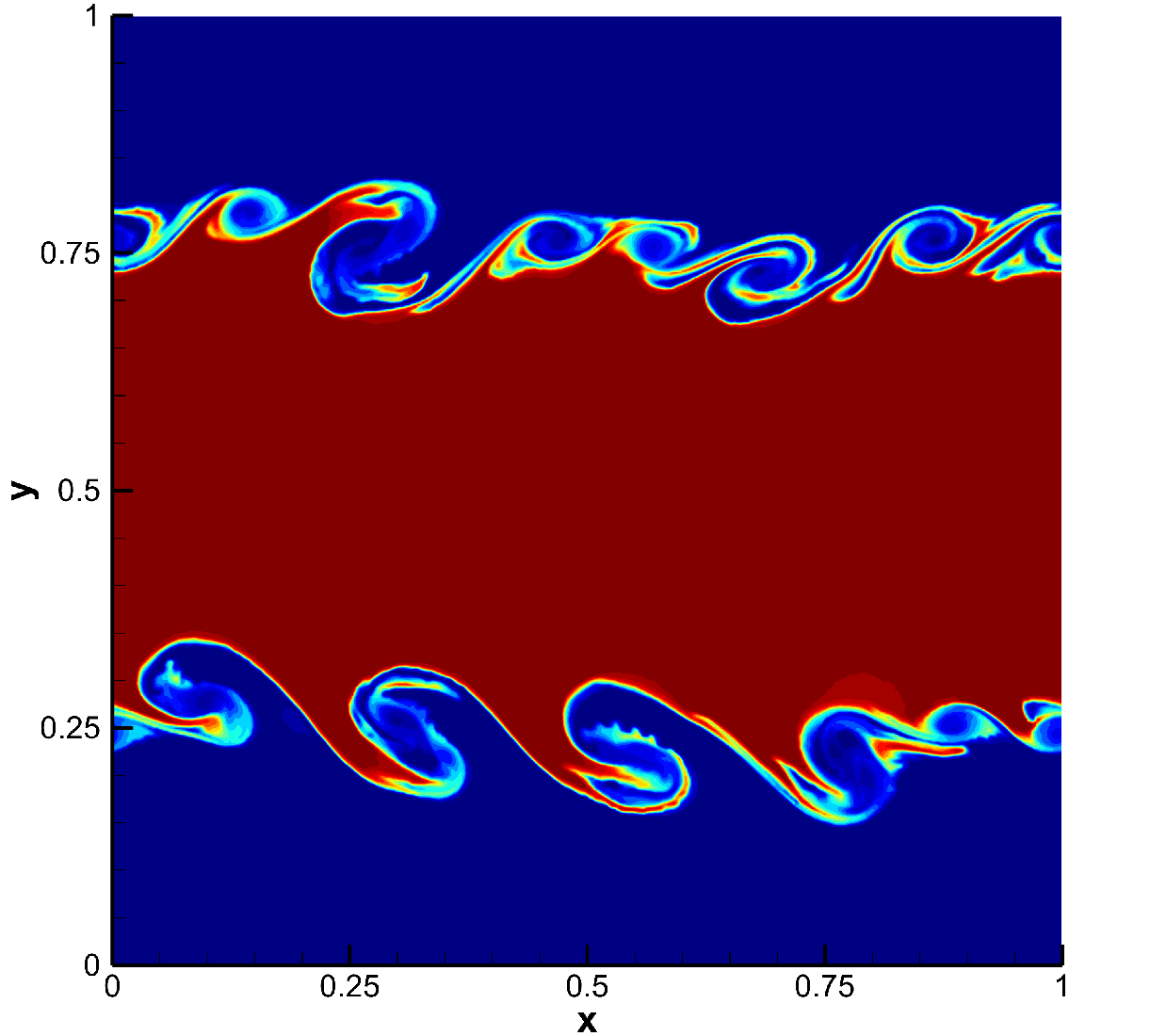}
\includegraphics[width=0.325\textwidth]{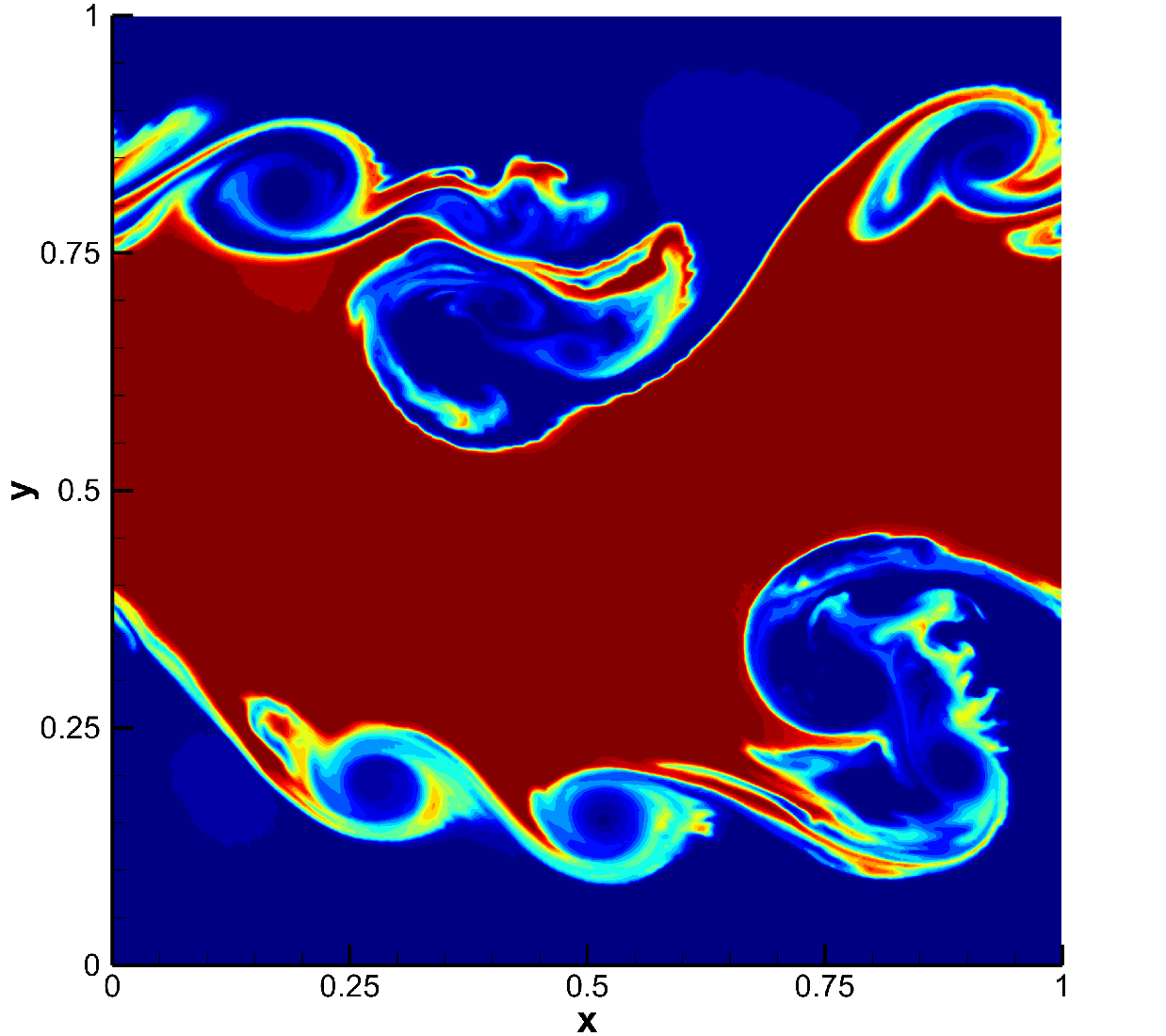}
\includegraphics[width=0.325\textwidth]{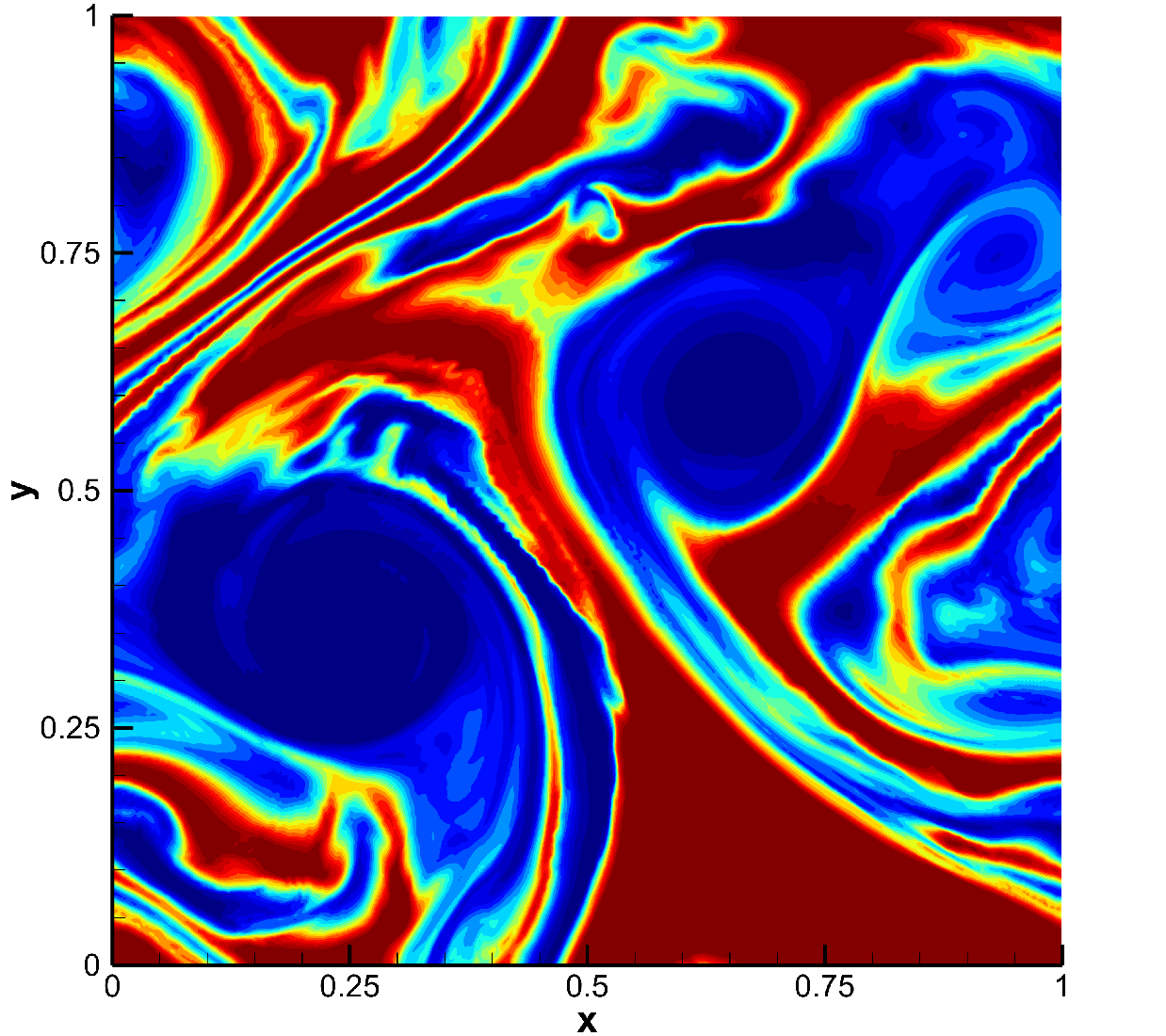}
\caption{\label{KH-instability-1} Kelvin--Helmholtz instability problem: Density contour plots obtained by the sixth-order FD-CGKS. Results on Grid 1 (top) and Grid 2 (bottom) are shown at \(t=1\), \(2.5\), and \(5\), from left to right. }
\end{figure}

The computation is performed using the present FD-CGKS(6-4) on two successively refined meshes: Grid~1 ($256\times 256$) and Grid~2 ($512\times 512$). Fig.~\ref{KH-instability-1} presents the density contours obtained on both grids at $t=1,~2$ and $5$. On the coarser Grid~1, the primary KH vortex roll-up and the secondary instabilities along the braids are clearly captured. As the mesh is refined to Grid~2, richer small-scale vortical structures are resolved, including finer secondary roll-ups as well as the thin filamentary features within the mixing layer.

\begin{figure}[!htb]
\centering
\includegraphics[width=0.65\textwidth]{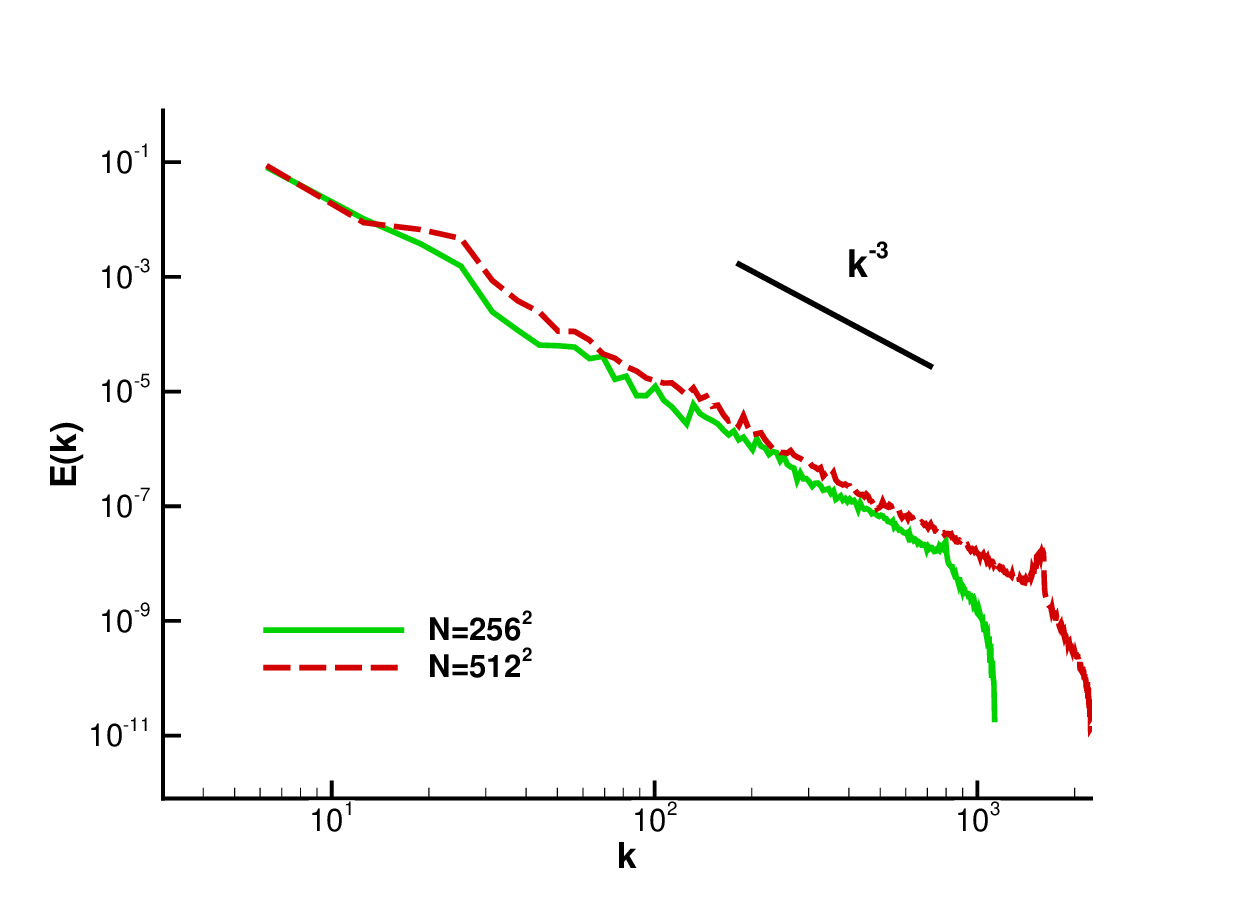}
\caption{\label{KH-instability-2} Kelvin--Helmholtz instability problem: Angle-averaged energy spectra at \( t=5 \) of different grid resolutions. Grid refinement broadens the correctly resolved spectrum, and both meshes accurately reproduce the \( k^{-3} \) inertial subrange scaling predicted by two-dimensional turbulence theory.}
\end{figure}

To further assess the resolution and numerical dissipation of the proposed scheme, its spectral behavior is examined.
Fig.~\ref{KH-instability-2} presents the angle-averaged kinetic energy spectra of the KH instability at $t = 5$, computed on meshes of different resolutions. Within the inertial subrange, the computed
spectra closely follow the classical $k^{-3}$ decay rate, in excellent agreement with the Kraichnan--Batchelor--Leith (KBL) theory for two-dimensional turbulence~\cite{san2015-KHI}. As the mesh is refined, this scaling extends to progressively higher wavenumbers, indicating that the effective numerical dissipation decreases and that an increasingly broad range of small-scale structures is resolved. These results confirm that the present FD-CGKS(6-4) features low numerical dissipation and high spectral resolution, enabling it to effectively capture the multi-scale flow structures.

\begin{figure}[!htb]
\centering
\includegraphics[width=0.495\textwidth]{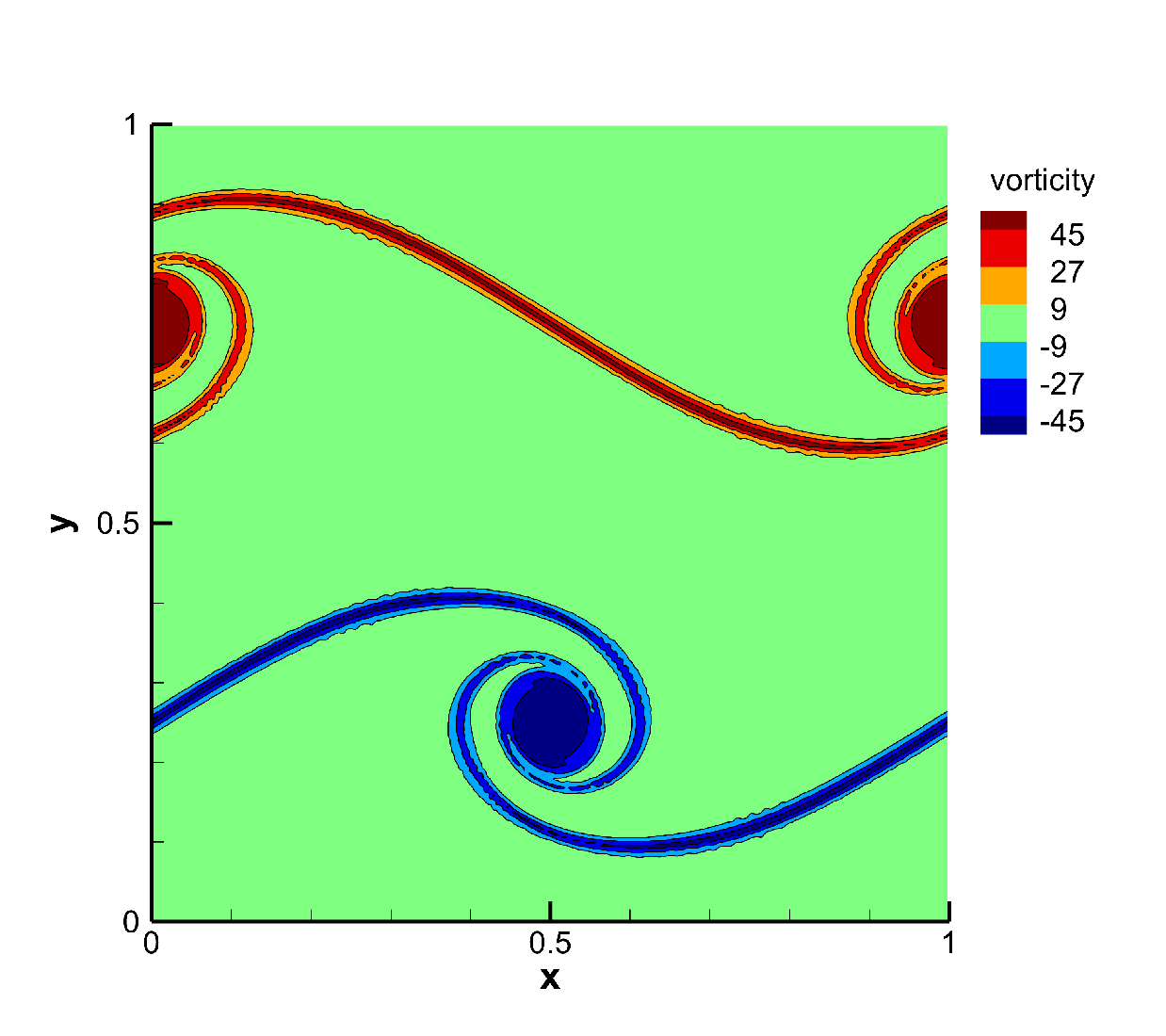}
\includegraphics[width=0.495\textwidth]{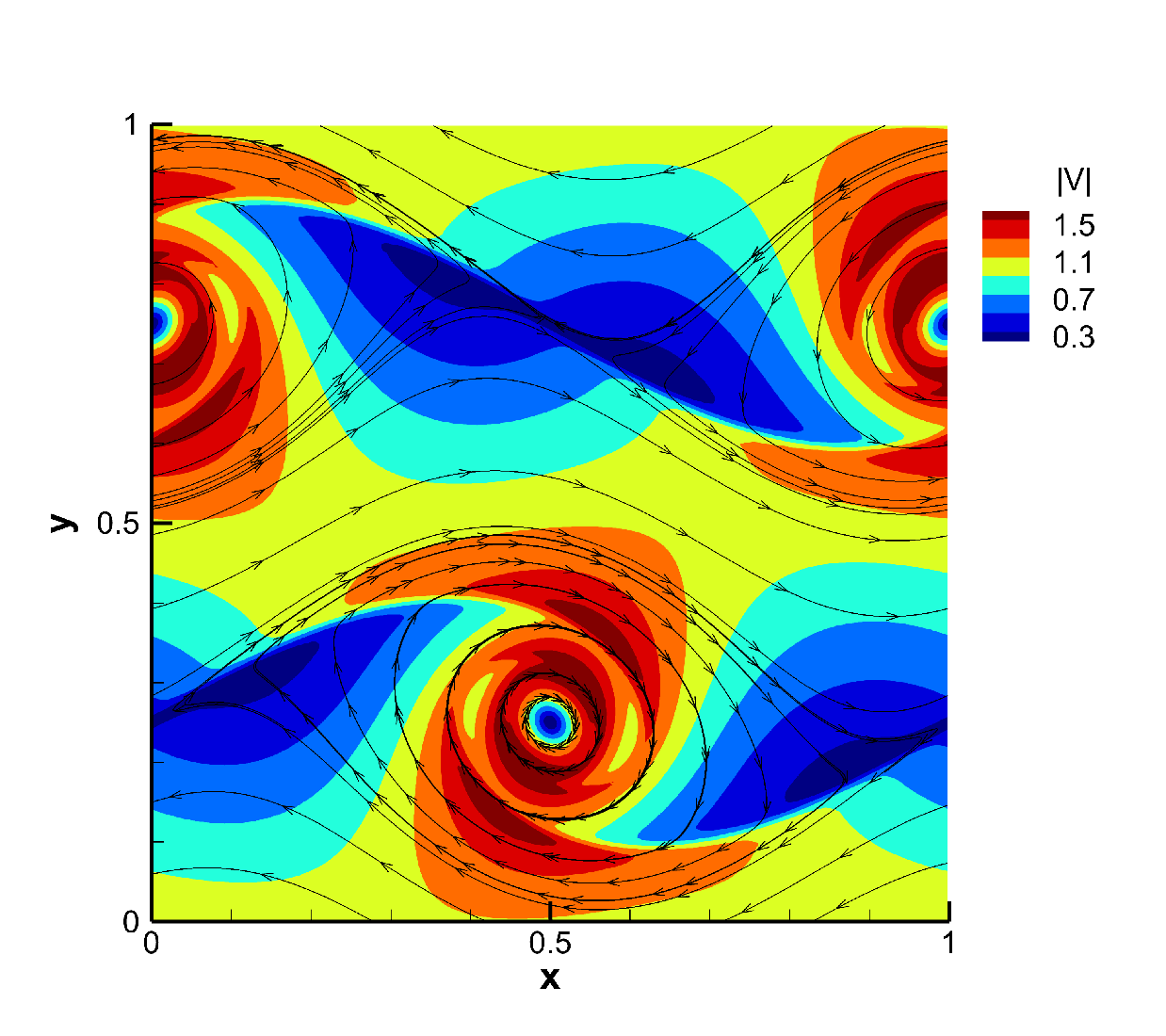}
\caption{\label{Double-shear-layer-1} Double shear layer problem: Vorticity contours (left) and velocity magnitude contours overlaid with streamlines (right), computed using the nonlinear sixth-order FD-CGKS. }
\end{figure}

\begin{figure}[!htb]
\centering
\includegraphics[width=0.495\textwidth]{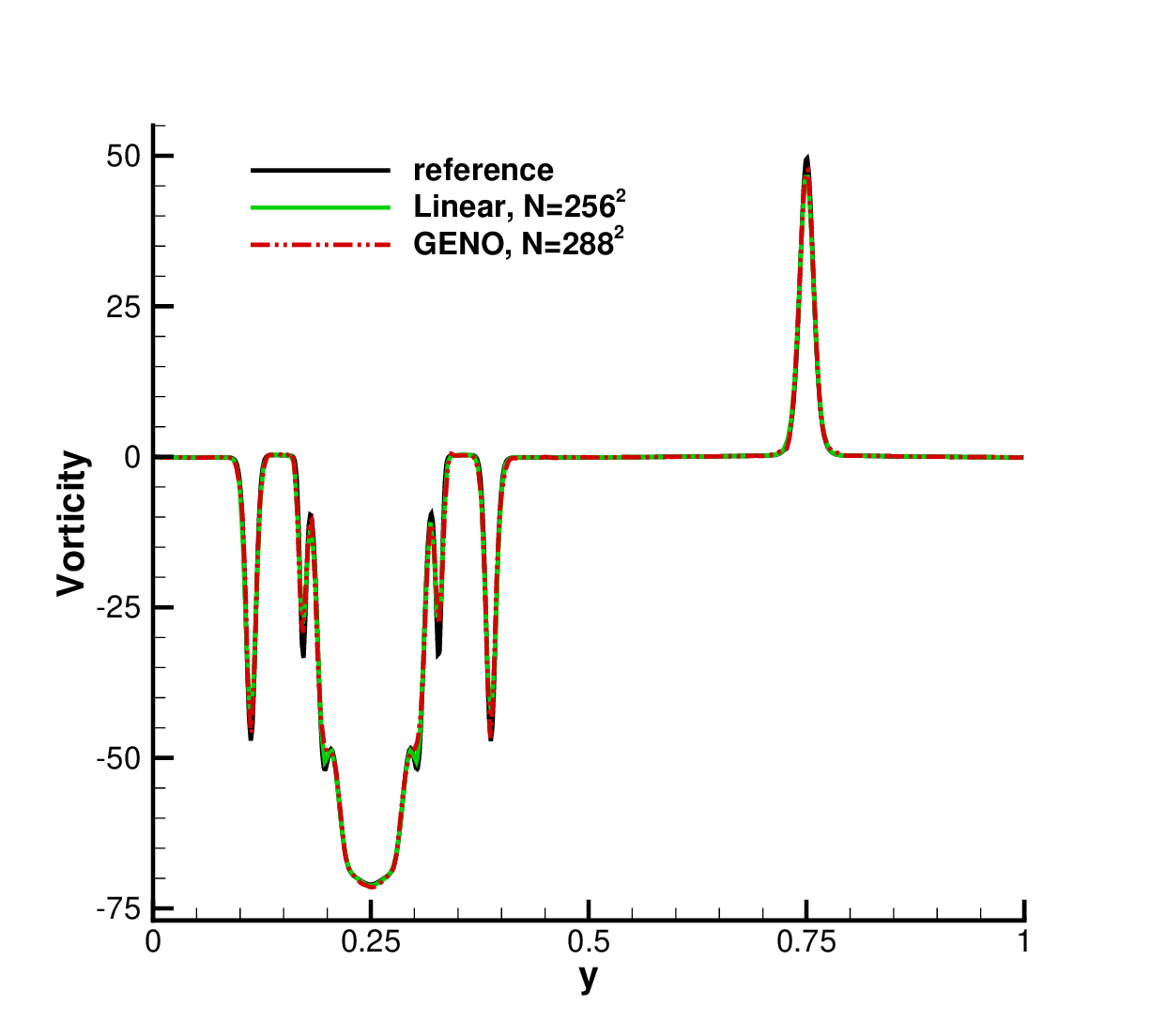}
\includegraphics[width=0.495\textwidth]{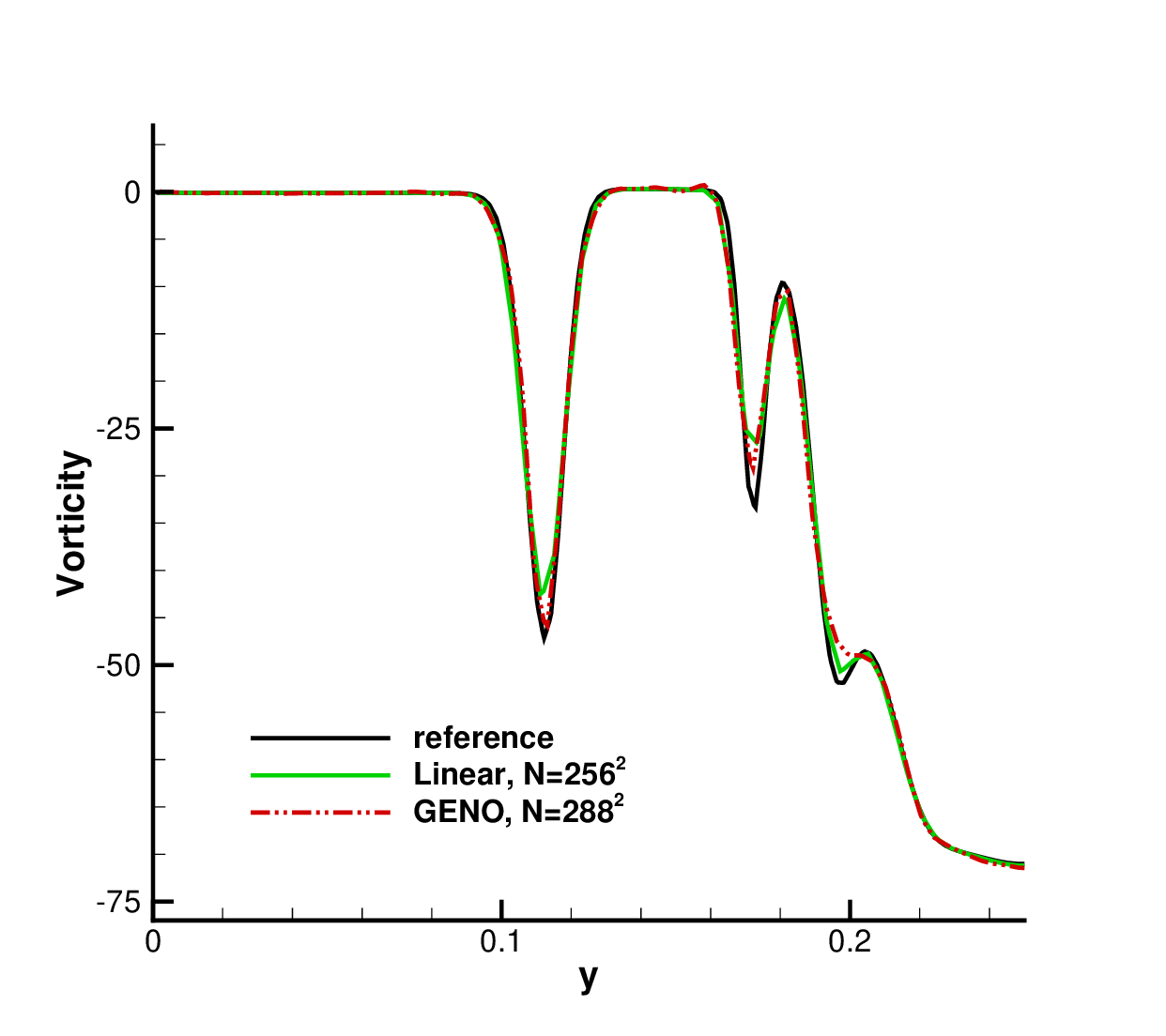}\\
\caption{\label{Double-shear-layer-2} Double shear layer problem: Vorticity profiles extracted along the centerline \( x=0.5 \) (left) and a close-up view (right), comparing the linear and nonlinear variants of the sixth-order FD-CGKS. }
\end{figure}

\subsection{Double shear layer problem}

The double shear layer is a canonical benchmark for assessing the accuracy and resolution of numerical schemes in the nearly incompressible regime. Since the thin shear layers are highly sensitive
to numerical perturbations and tend to develop spurious vortices when under-resolved, this problem provides a stringent test of both the accuracy and the stability of the proposed scheme.
Following \cite{Brown1995}, the initial condition is prescribed as
\begin{align*}
  U(x,y)=
  \begin{cases}
    \tanh\!\bigl(k(y-0.25)\bigr), & y\le 0.5,\\[4pt]
    \tanh\!\bigl(k(0.75-y)\bigr), & y>0.5,
  \end{cases}
  \qquad
  V(x,y)=\delta\sin(2\pi x),
\end{align*}
with uniform density and pressure
\begin{align*}
  \rho(x,y)=1, \qquad
  p(x,y)=\frac{\rho}{M_a^{2}\,\gamma},
\end{align*}
where $k=100$, $\delta=0.05$, and $M_a=0.15$. The kinematic viscosity is $\nu=5.0\times10^{-5}$.
The computational domain is $[0,1]\times[0,1]$ with periodic boundary conditions in both the $x$- and $y$-directions. The final computation time is set to $t=0.8$.

Numerical results are obtained using both the linear and nonlinear variants of the sixth-order FD-CGKS(6-4), where the respective linear or nonlinear strategies are consistently applied to both the spatial reconstruction and the finite-difference approximation of the flux derivatives.
The computations for the linear and nonlinear schemes are performed on uniform meshes of \( 256 \times 256 \) and \( 288 \times 288 \) grid points, respectively. Fig.~\ref{Double-shear-layer-1} presents the flow fields obtained by the nonlinear FD-CGKS(6-4). The left panel displays the vorticity contours, while the right panel shows streamlines overlaid on the velocity magnitude contours.
The computed overall flow topology faithfully reproduces the reference solution \cite{zhao2019compact}, validating the scheme's low numerical dissipation and its stability against spurious solutions.

A quantitative evaluation is provided in Fig.~\ref{Double-shear-layer-2}, which compares the vorticity distribution along the centerline \( x=0.5 \) against the reference solution from \cite{zhao2019compact}. The right panel provides a close-up view. It is observed that the nonlinear FD-CGKS(6-4) introduces slightly larger deviations compared to its linear counterpart. This behavior is expected and is attributable to the additional numerical dissipation inherent to the nonlinear adaptation mechanism.

\begin{figure}[!htb]
\centering
\includegraphics[width=0.495\textwidth]{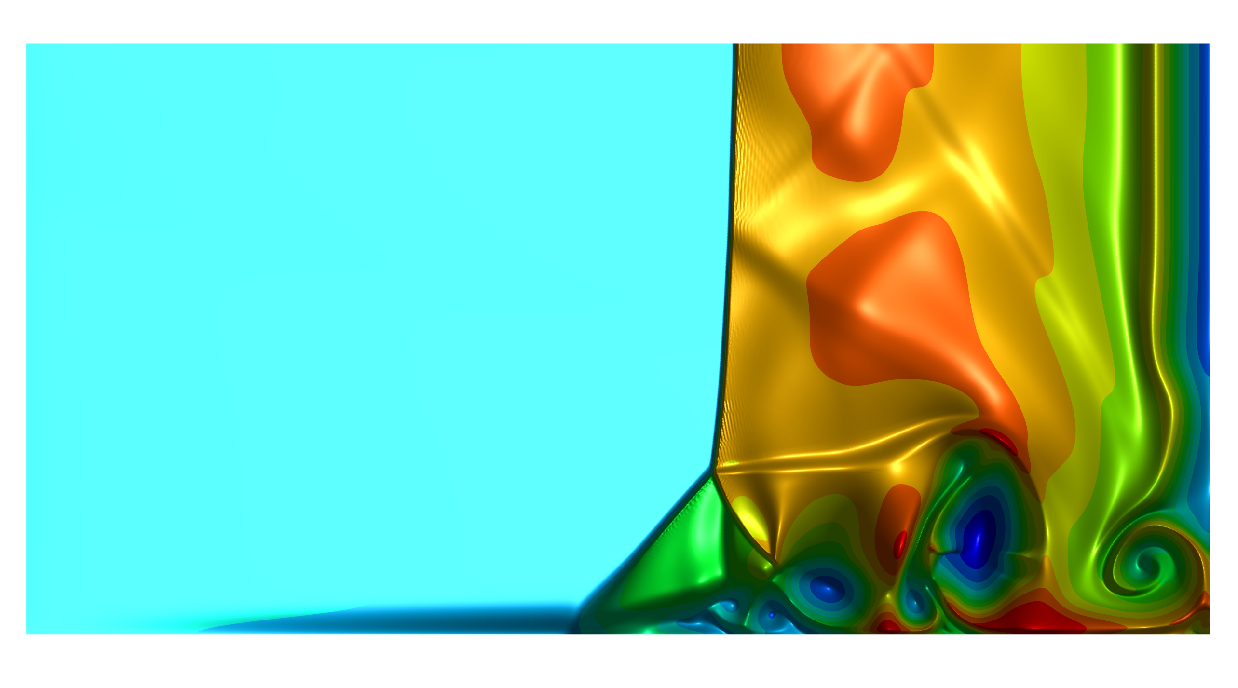}
\includegraphics[width=0.495\textwidth]{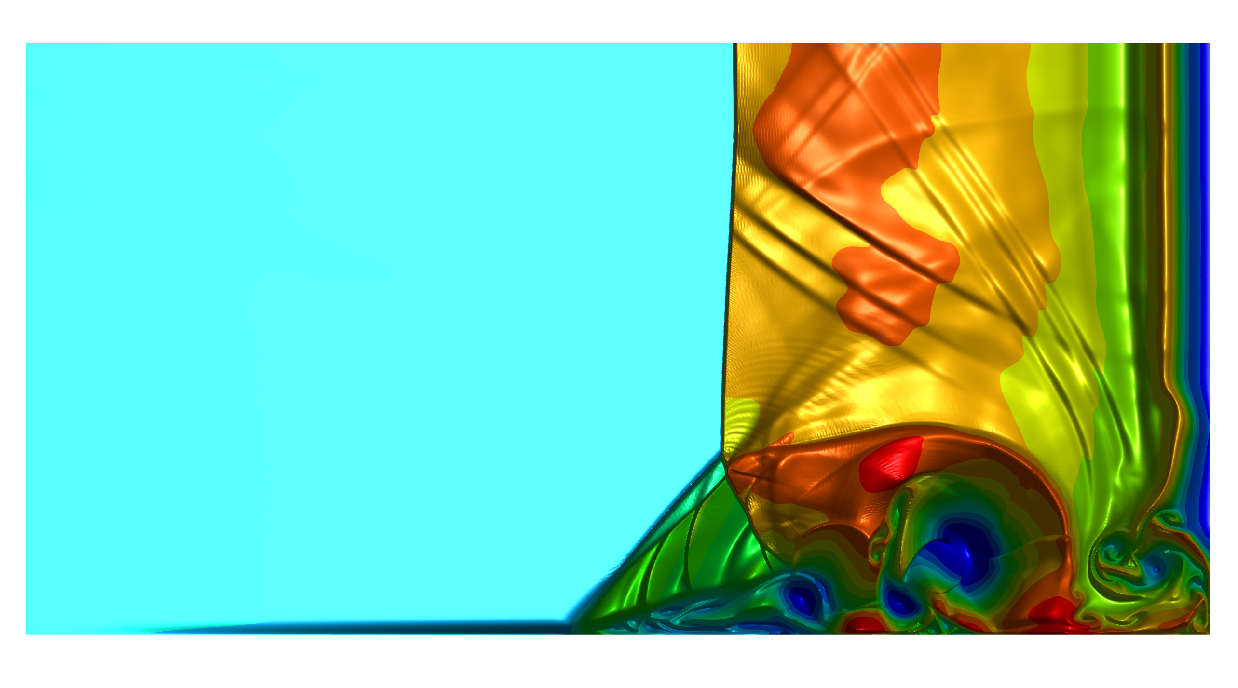}\\
\caption{\label{vis-shock-tube-1} Viscous shock tube problem: Density contours at \( t = 1.0 \) for \( \text{Re} = 200 \) (left) and \( \text{Re} = 1000 \) (right). Computed on a \( 2000 \times 1000 \) mesh. Eleven equally spaced contour levels from \( 25 \) to \( 110 \) are displayed. }
\end{figure}

\begin{figure}[!htb]
\centering
\includegraphics[width=0.495\textwidth]{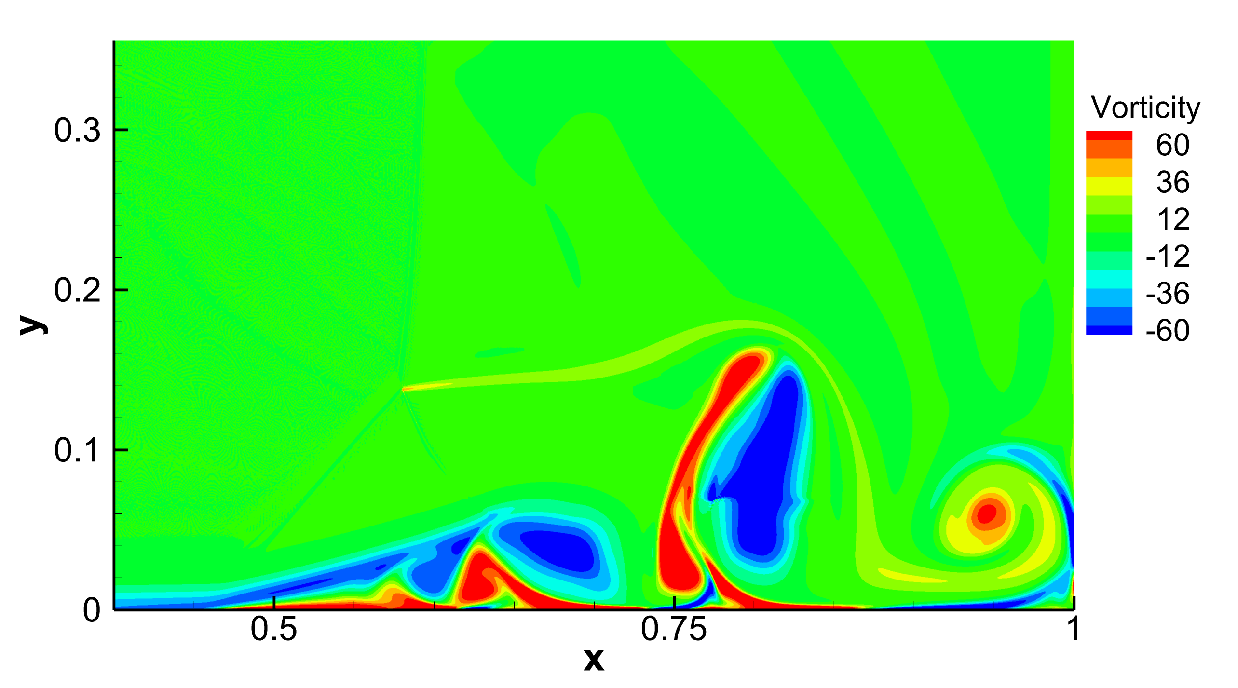}
\includegraphics[width=0.495\textwidth]{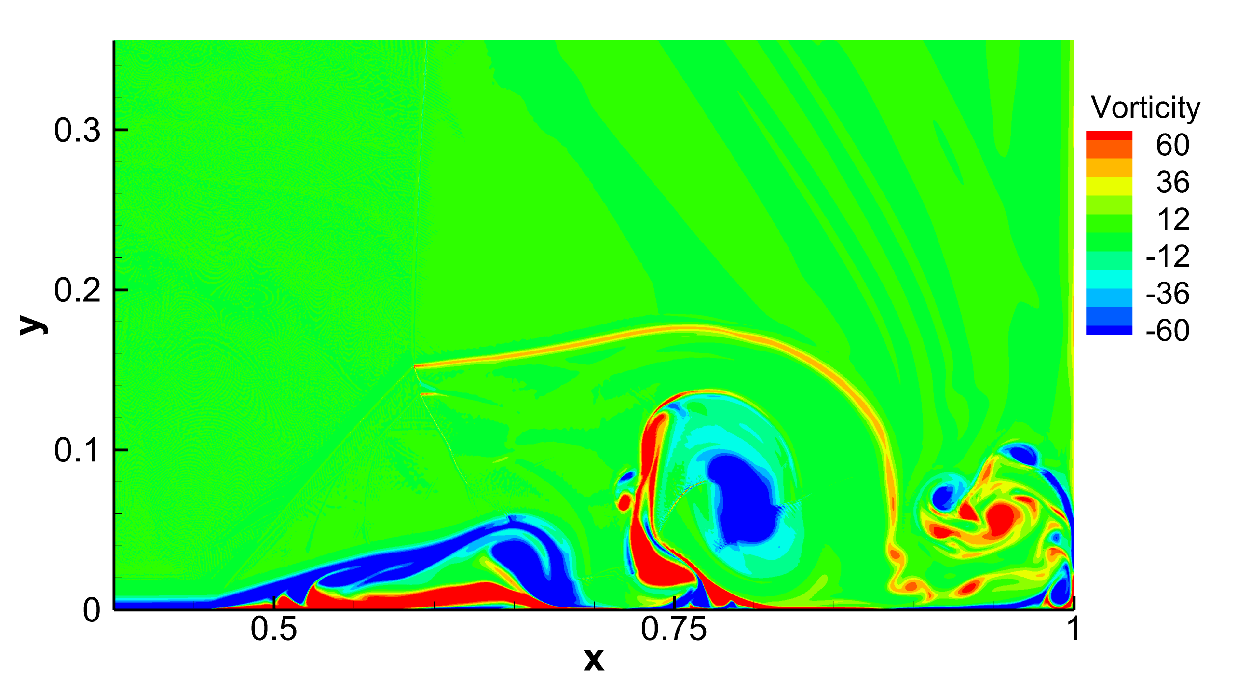}\\
\caption{\label{vis-shock-tube-2} Viscous shock tube problem: Vorticity contours at \( t = 1.0 \) for \( \text{Re} = 200 \) (left) and \( \text{Re} = 1000 \) (right). Computed on a \( 2000 \times 1000 \) mesh. }
\end{figure}

\subsection{Viscous shock tube problem}

The viscous shock tube problem features complex flow phenomena driven by intense shock-wave and boundary-layer interactions \cite{daru2000_test}. It serves as a rigorous benchmark to evaluate both the robustness and the high-resolution capability of numerical schemes in capturing intricate viscous structures. The flow is enclosed within a unit square cavity. Taking advantage of the geometric symmetry, the computational domain is restricted to the lower half, defined as \( [0, 1] \times [0, 0.5] \). A symmetric boundary condition is applied at the top boundary (\( y = 0.5 \)), while adiabatic no-slip wall conditions are imposed on the remaining boundaries (left, right, and bottom). The initial condition is given by
\begin{align*}
(\rho, U, V, p) =
\begin{cases}
(120, 0, 0, 120/\gamma), & 0 \leq x < 0.5, \\
(1.2, 0, 0, 1.2/\gamma), & 0.5 \leq x \leq 1.
\end{cases}
\end{align*}
The Prandtl number is set to \( \text{Pr} = 0.72 \), and the simulations are advanced to a final time of \( t = 1.0 \).

Two cases, corresponding to Reynolds numbers of \( \text{Re} = 200 \) and \( 1000 \) (with dynamic viscosities \( \mu = 0.005 \) and \( 0.0001 \), respectively), are simulated. For the low Reynolds number case (\( \text{Re} = 200 \)), uniform meshes of \( 500 \times 250 \), \( 1000 \times 500 \), and \( 2000 \times 1000 \) points are employed. For the high Reynolds number case (\( \text{Re} = 1000 \)), the meshes are refined to \( 1200 \times 600 \), \( 1600 \times 800 \), and \( 2000 \times 1000 \) points. Reference solutions are taken from \cite{zhou2018grid}, which were computed using a fifth-order finite volume scheme on highly refined grids.
Figs.~\ref{vis-shock-tube-1} and \ref{vis-shock-tube-2} display the density and vorticity contours for \( \text{Re} = 200 \) (left column) and \( \text{Re} = 1000 \) (right column), respectively, computed on the present finest mesh of \( 2000 \times 1000 \) points. The proposed FD-CGKS(6-4) successfully resolves the rich vortex structures generated by the interaction, without exhibiting spurious numerical oscillations near the shock waves.

\begin{figure}[!htb]
\centering
\includegraphics[width=0.495\textwidth]{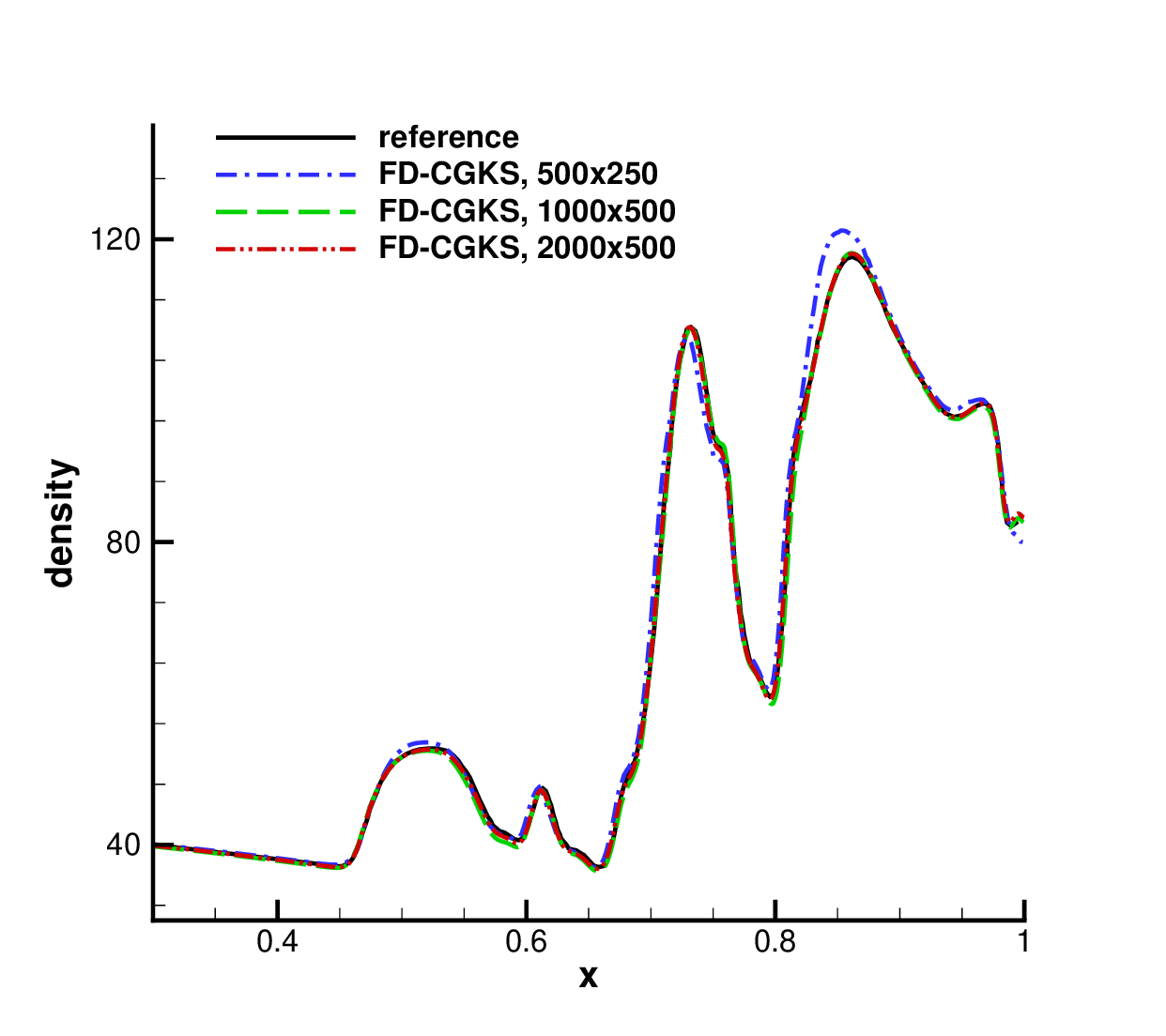}
\includegraphics[width=0.495\textwidth]{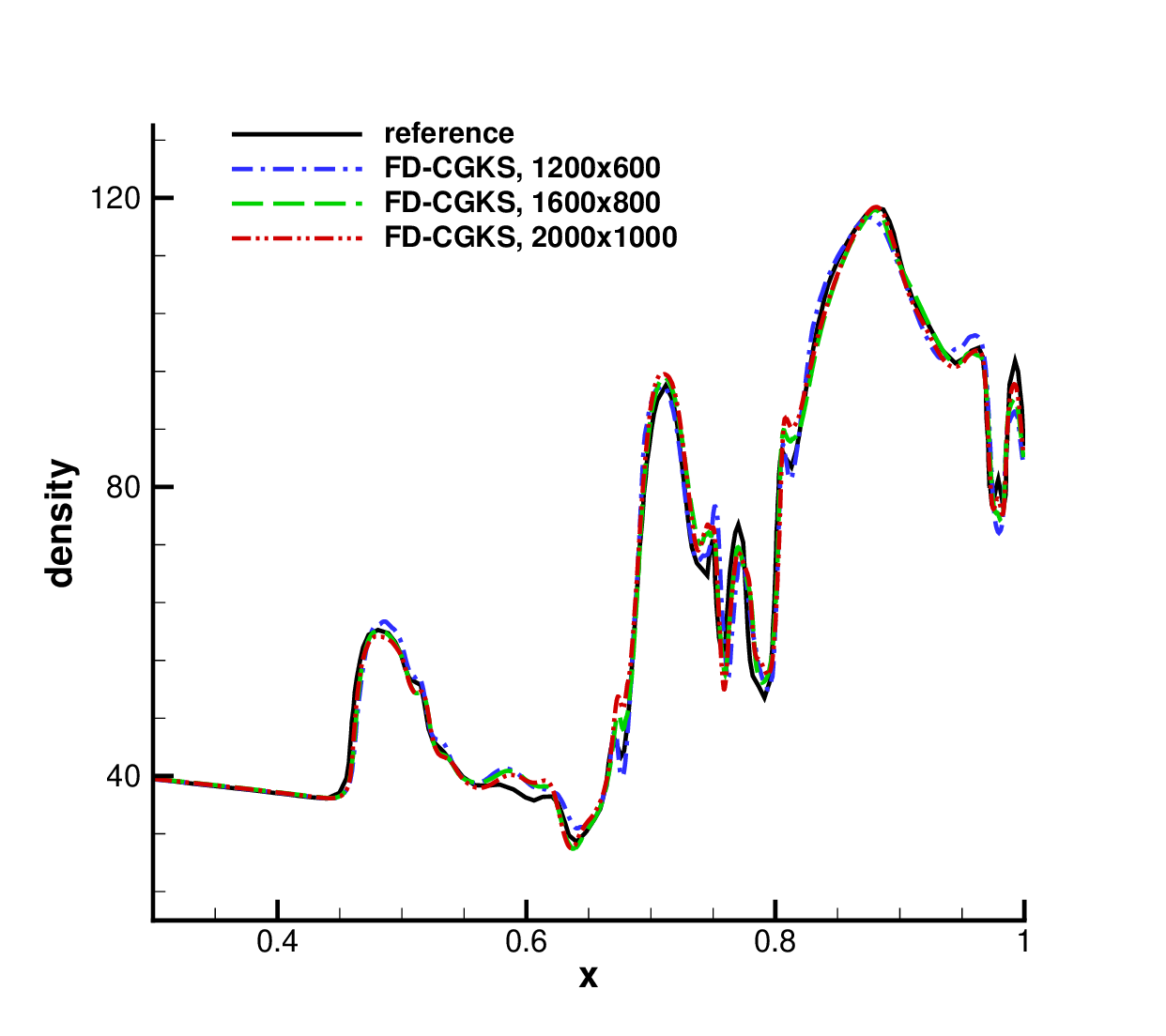}
\caption{\label{vis-shock-tube-3} Viscous shock tube problem: Density profiles along the bottom wall (\( y = 0 \)) at \( t = 1.0 \) for \( \text{Re} = 200 \) (left) and \( \text{Re} = 1000 \) (right). }
\end{figure}

\begin{figure}[!htb]
\centering
\includegraphics[width=0.495\textwidth]{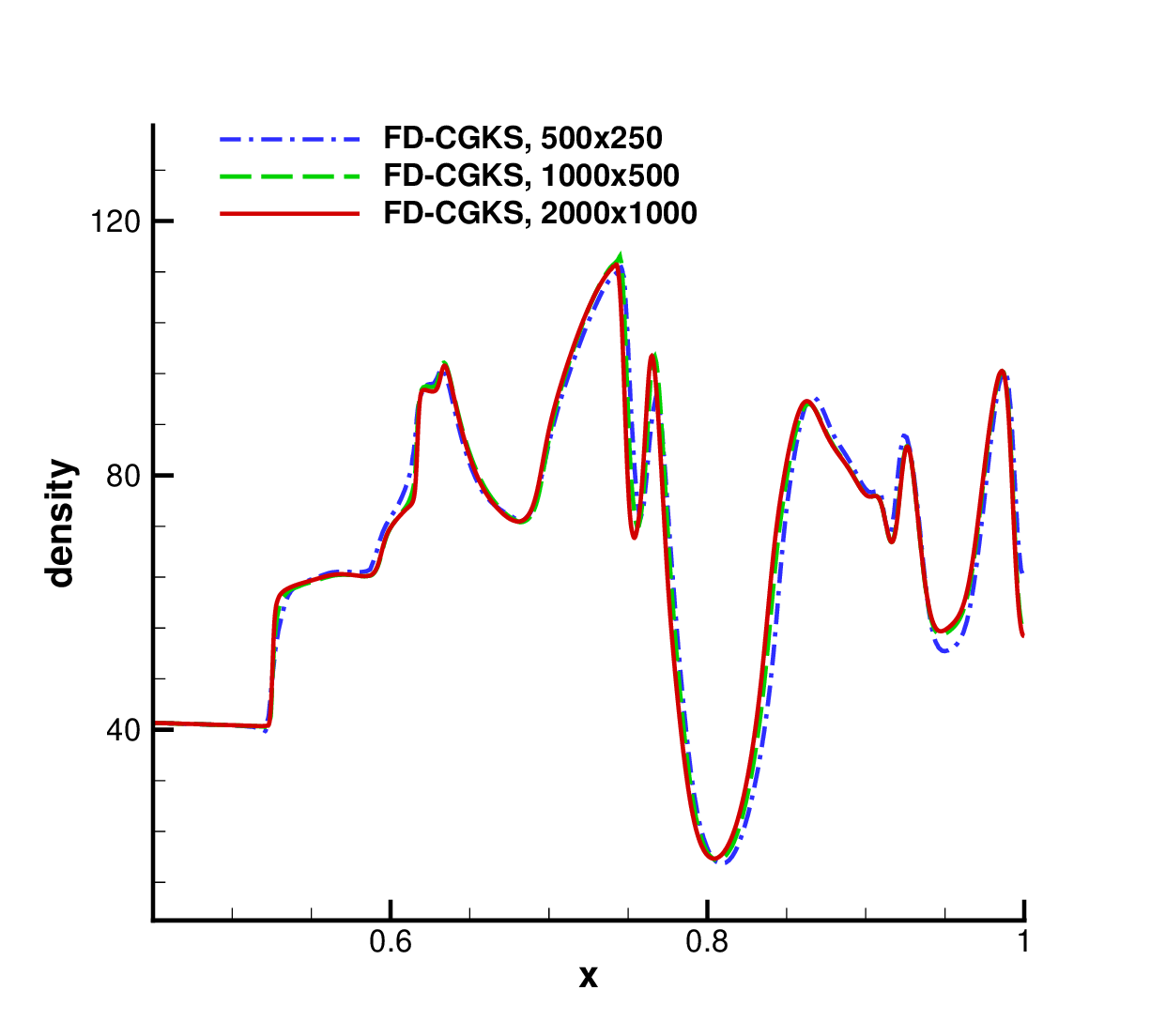}
\includegraphics[width=0.495\textwidth]{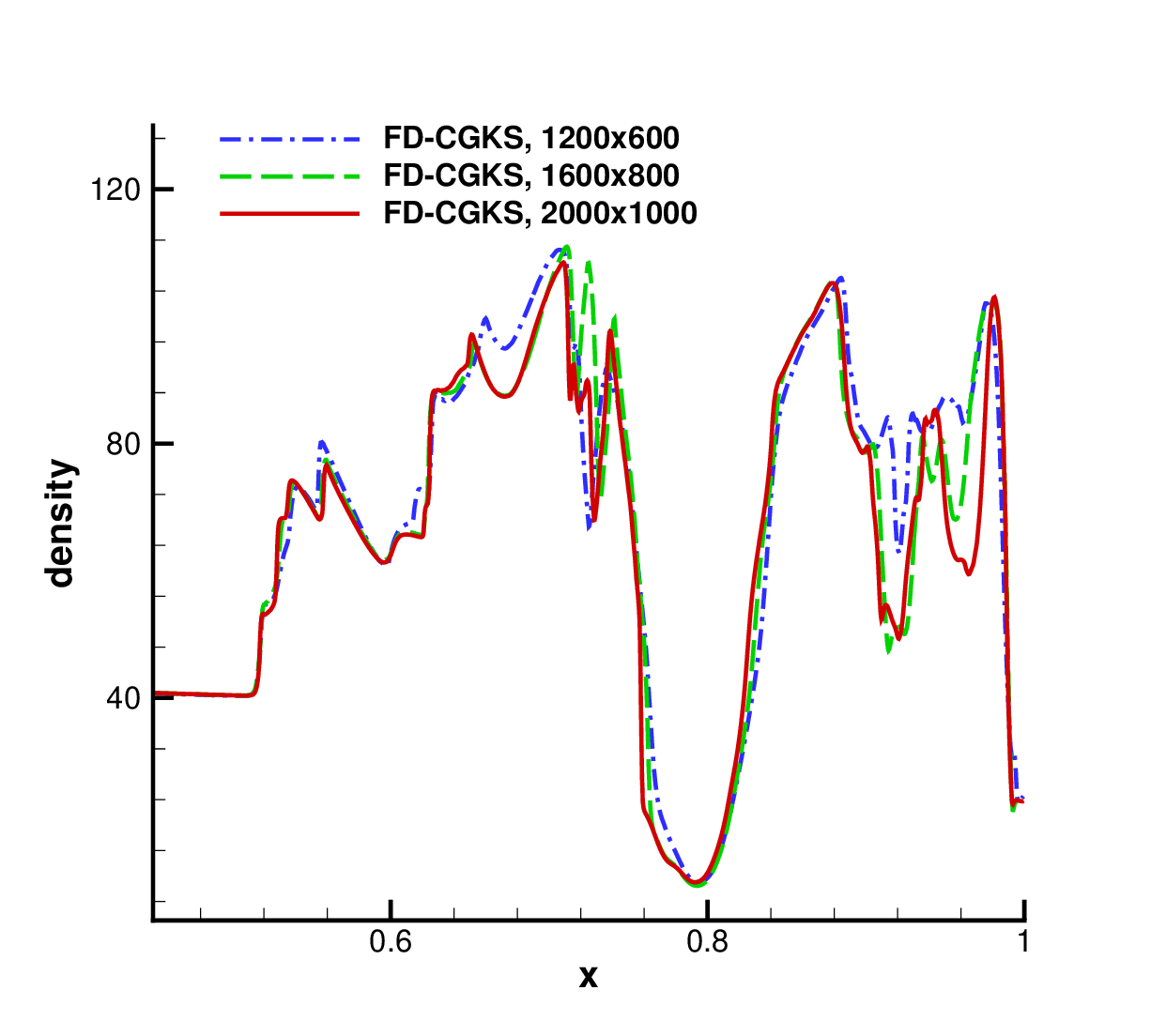}
\caption{\label{vis-shock-tube-4} Viscous shock tube problem: Density profiles along the horizontal line \( y = 0.075 \) at \( t = 1.0 \) for \( \text{Re} = 200 \) (left) and \( \text{Re} = 1000 \) (right). }
\end{figure}

Figs.~\ref{vis-shock-tube-3} and \ref{vis-shock-tube-4} provide quantitative comparisons of the density distributions extracted along the bottom wall and the horizontal line \( y=0.075 \), respectively.
For \( \text{Re} = 200 \), the coarsest grid already achieves satisfactory grid convergence, with only slight deviations near the local extrema of the density profile. The results on the finer meshes align closely with both each other and the reference solution.
For \( \text{Re} = 1000 \), the density profiles along the bottom wall computed by the FD-CGKS(6-4) demonstrate overall consistency across the three meshes. As the grid is refined, the numerical solutions asymptotically approach the reference data, with discrepancies confined to the highly localized, smaller-scale structures. These minor deviations are attributable to the still-insufficient grid resolution; the number of grid points in the present finest mesh is only about \( 1/6.25 \) of that used for the reference solution.
Along the line \( y=0.075 \), the density distributions exhibit more pronounced variations among the three grids, reflecting the presence of even finer-scale flow structures in this region. This stringent test case further substantiates the numerical robustness and high resolving power of the proposed FD-CGKS for simulating viscous flows characterized by strong shock-wave and boundary-layer interactions.

\section{Conclusion}
\label{sec:conclusion}

In this paper, a high-order FD-CGKS has been developed and systematically validated for the simulation of compressible inviscid and viscous flows. A primary novelty of this work is the introduction of a dual-mesh-based finite difference framework for structured meshes. Unlike traditional finite difference methods, the proposed approach utilizes conservative variables at both the main and dual nodes. This design significantly simplifies the multidimensional spatial reconstruction required for the space-time coupled flux evaluations inherent to the GKS. By exploiting the time-accurate interface solutions provided by the gas-kinetic model, the scheme explicitly updates the averaged derivatives, thereby enabling highly efficient and compact high-order reconstruction. Furthermore, to ensure robustness and maintain non-oscillatory properties near strong shocks, a nonlinear GENO method has been successfully integrated into the spatial reconstruction and numerical flux evaluation processes.

The performance of the proposed FD-CGKS was rigorously evaluated through a variety of one- and two-dimensional benchmark problems. Theoretical accuracy tests confirmed the scheme's expected high-order convergence. Simulations of challenging one-dimensional problems, such as interacting blast waves and shock-density wave interactions, demonstrated its robust shock-capturing capabilities and high spatial resolution. In comparative studies among one-dimensional schemes, the FD-CGKS(6-4) variant exhibited a superior balance between numerical accuracy and computational robustness. In two-dimensional viscous and inviscid tests, including the double Mach reflection, Kelvin-Helmholtz (KH) instability, and viscous shock tube problems, the scheme exhibited excellent low-dissipation characteristics.
Notably, spectral analysis of the KH instability demonstrated that, upon grid refinement, the scheme accurately recovers the classical $k^{-3}$ scaling law of two-dimensional turbulence in the inertial subrange. Additionally, the viscous shock tube results highlighted the scheme's ability to accurately resolve intricate shock-wave/boundary-layer interactions without spurious oscillations, achieving high-fidelity results on significantly coarser meshes than those required by traditional methods.

In summary, the proposed FD-CGKS achieves an outstanding balance between numerical accuracy, low dissipation, and shock-capturing robustness. This makes it a highly promising numerical tool for simulating complex compressible flows characterized by multiscale structures. Future work will focus on extending this dual-mesh finite difference framework to large-scale simulations of compressible turbulence.

\section*{Acknowledgments}
The current research is supported by National Key R\&D Program of China (Grant Nos. 2022YFA1004500), National Science Foundation of China (92371107), and Hong Kong research grant council (16208324).

\appendix
\section{Gas evolution model and evaluation of flow variables and fluxes}
\label{appendix:A}

\subsection{Time-accurate solution of the gas distribution function}
\label{appendix:A1}

In the GKS, the evolution of the gas distribution function is governed by the kinetic BGK equation,
\begin{equation*}
    f_t + \mathbf{u} \cdot \nabla f = \frac{g - f}{\tau},
\end{equation*}
where \(\mathbf{u}\) is the particle velocity, \(f\) is the gas distribution function, \(g\) is the corresponding equilibrium state (a Maxwellian distribution) that \(f\) approaches, and \(\tau\) is the particle collision time.
To evaluate the flow variables and fluxes at the grid interface, the time-accurate gas distribution function \(f\) is obtained.
Based on the BGK equation, its integral solution provides an implicit form for \(f\) along the particle trajectory:
\begin{equation}\label{eq:integral_f}
    f(\mathbf{x}_0, t, \mathbf{u}, \xi) = \frac{1}{\tau} \int_0^t g(\mathbf{x}', t', \mathbf{u}, \xi) e^{-(t-t')/\tau} \mathrm{d}t' + e^{-t/\tau} f_0(\mathbf{x}_0 - \mathbf{u}(t-t_0), \mathbf{u}, \xi),
\end{equation}
where \(\mathbf{x}_0\) is the numerical quadrature point at the interface, \(\mathbf{x}_0 = \mathbf{x}' + \mathbf{u}(t - t')\) is the particle trajectory, and \(f_0\) is the initial state of the gas distribution function at \(t = 0\). By modeling \(f_0\) and \(g\), Eq.~\eqref{eq:integral_f} forms the basis for constructing the time-dependent distribution function at the cell interface, and the solution $f(\mathbf{x}_0, t, \mathbf{u}, \xi)$ is fully determined by the macroscopic flow state at interfaces \cite{xu2}.

To connect the gas distribution function to the macroscopic variables, we first define the collision invariant vector $\boldsymbol{\psi}$ and the integration measure in the phase space $\mathrm{d}\Xi$. For a 2D flow, they are given by
\begin{align*}
\begin{split}
    &\boldsymbol{\psi} = (\psi_1, \psi_2, \psi_3, \psi_4)^T = \left(1, u, v, \frac{1}{2}(u^2 + v^2 + \xi^2)\right)^T,\\
    &\mathrm{d}\Xi = \mathrm{d}u \mathrm{d}v \mathrm{d}\xi_1 ... \mathrm{d}\xi_K,
\end{split}
\end{align*}
where $u$ and $v$ are the particle velocities in the $x$- and $y$-directions, respectively. $\xi$ is the random velocity with $\xi^2 = \xi_1^2 + \xi_2^2 + ... + \xi_K^2$, and $K$ is the number of internal degrees of freedom.
The macroscopic conservative variables, including mass $\rho$, momentum $(\rho U, \rho V)$, and energy $\rho E$, can be evaluated by taking the moments of the gas distribution function $f$:
\begin{equation}
    \mathbf{W} =
    \begin{pmatrix}
        \rho \\
        \rho U \\
        \rho V \\
        \rho E
    \end{pmatrix}
    = \int f \boldsymbol{\psi} \mathrm{d}\Xi.
\end{equation}
Similarly, the corresponding fluxes for mass, momentum, and energy across a grid interface in the normal direction (aligned with velocity $u$) can be obtained by taking the moments of $u f$:
\begin{equation}
    \mathbf{F} = \int u f \boldsymbol{\psi} \mathrm{d}\Xi.
\end{equation}

In the numerical implementation, especially for the second-order evolution model, the time-accurate distribution function $f$ at the grid interface can be approximated as a linear function of time:
\begin{equation}\label{eq:f_linear}
    \hat{f}(t) = f^n + t f_t^n.
\end{equation}
By defining the time integration of $f(t)$ as $\bar{f}(t) = \int_0^t f(t') \mathrm{d}t'$, the unknowns $f^n$ and $f_t^n$ in Eq.~\eqref{eq:f_linear} can be determined using the integrated distribution function at different time stages (e.g., $\Delta t/2$ and $\Delta t$):
\begin{equation}
\begin{aligned}
    f^n &= \frac{1}{\Delta t} \left( 4\bar{f}(\Delta t/2) - \bar{f}(\Delta t) \right), \\
    f_t^n &= \frac{4}{\Delta t^2} \left( -2\bar{f}(\Delta t/2) + \bar{f}(\Delta t) \right).
\end{aligned}
\end{equation}

\subsection{Time-accurate solution of the macroscopic variables}
\label{appendix:A2}

By taking the moments of \(\hat{f}(t)\) and its time derivative with respect to \(u\boldsymbol{\psi}\) at \(t=0\), the numerical flux and its time derivative used for the macroscopic variable updates are evaluated as:
\begin{equation}
\begin{aligned}
    \widehat{\mathbf{F}} &= \int u f^n \boldsymbol{\psi} \mathrm{d}\Xi, \\
    \widehat{\mathbf{F}}_t &= \int u f_t^n \boldsymbol{\psi} \mathrm{d}\Xi.
\end{aligned}
\end{equation}
For the sake of brevity, the time superscript \(n\) (denoting the time level \(t^n\)) for the macroscopic variables is omitted hereafter.

By taking the moments of \(\hat{f}(t)\) and its time derivative with respect to \(\boldsymbol{\psi}\) at \(t=0\), the macroscopic flow variables and their time derivatives are evaluated as:
\begin{align*}
\begin{split}
    \mathbf{W}^e &= \int f^n \boldsymbol{\psi} \mathrm{d}\Xi, \\
    \mathbf{W}^e_t &= \int f_t^n \boldsymbol{\psi} \mathrm{d}\Xi.
\end{split}
\end{align*}
Based on the evolution model presented in \cite{zhao2023direct}, the final evaluated interface flow variables used for the averaged gradient updates in the 1D case are given by:
\begin{equation}\label{eq:W-model}
\begin{split}
    \widehat{\mathbf{W}}^{l,r}   &= (1 - e^{-\Delta t/\tau_0})\mathbf{W}^e   + e^{-\Delta t/\tau_0}\mathbf{W}^{l,r}, \\
    \widehat{\mathbf{W}}^{l,r}_t &= (1 - e^{-\Delta t/\tau_0})\mathbf{W}^e_t + e^{-\Delta t/\tau_0}\mathbf{W}^{l,r}_t,
\end{split}
\end{equation}
where \(\mathbf{W}^{l,r}\) are evaluated from the independently evolved non-equilibrium gas distribution functions, \(f^{l,r}(t)\), on the left and right sides of the grid interface (see \cite{zhao2023direct} for further details).

\subsection{Two-dimensional solutions of flow variables}
\label{appendix:A3}

For the 2D case, two pairs of discontinuous left and right states are obtained based on the initial distributions along the two respective directions, denoted as \(\widehat{\mathbf{W}}^{l,(1)}\), \(\widehat{\mathbf{W}}^{r,(1)}\) and \(\widehat{\mathbf{W}}^{l,(2)}\), \(\widehat{\mathbf{W}}^{r,(2)}\). For a grid interface whose normal vector aligns with direction \((1)\), the final evaluated solutions at the interface, which are used for the averaged gradient updates, are constructed as follows:
\begin{equation}\label{eq:W-model-2d}
\begin{split}
    \widetilde{\mathbf{W}}^{l} &= \frac{1}{2}\widehat{\mathbf{W}}^{l,(1)} + \frac{1}{4}\left(\widehat{\mathbf{W}}^{l,(2)} + \widehat{\mathbf{W}}^{r,(2)}\right), \\
    \widetilde{\mathbf{W}}^{r} &= \frac{1}{2}\widehat{\mathbf{W}}^{r,(1)} + \frac{1}{4}\left(\widehat{\mathbf{W}}^{l,(2)} + \widehat{\mathbf{W}}^{r,(2)}\right).
\end{split}
\end{equation}
The corresponding time derivatives, \(\widetilde{\mathbf{W}}^{l}_t\) and \(\widetilde{\mathbf{W}}^{r}_t\), are evaluated in an analogous manner.
For a grid interface whose normal vector aligns with direction \((2)\), similar expressions can be readily deduced.
Furthermore, Eq. (\ref{eq:W-model-2d}) is specifically constructed to better preserve the multidimensionality of the evolved solutions.

\section{Sub-stencil polynomials and smoothness indicators}
\label{appendix:C}

In this section, as the basic lower-order part of the GENO reconstruction, we provide the detailed expressions for the sub-stencil polynomials and their corresponding smoothness indicators used in the compact reconstruction at both the grid interfaces and the solution nodes.
Here, \(h= \Delta x\) denotes the mesh size, \(Q_k\) and \(Q_{x,k}\) represent the macroscopic variables and their spatially averaged derivatives over the interval $[x_{k-1/2},x_{k+1/2}]$, respectively.

\subsection{Reconstruction at interfaces (left side)}
\label{appendix:C1}

For the left-side reconstruction at the grid interface \(j+1/2\), the sub-stencil polynomials for the value \(q_{k}^l\) and the derivative \(q_{k,x}^l\) are given by:
\begin{equation}
\begin{aligned}
    q_{1}^l &= \frac{1}{4} \left( -5Q_{j-1} + 9Q_j - 3Q_{x,j-1}h \right), \\
    q_{2}^l &= \frac{1}{8} \left( -Q_{j-1} + 6Q_j + 3Q_{j+1} \right), \\
    q_{3}^l &= \frac{1}{8} \left( 3Q_j + 6Q_{j+1} - Q_{j+2} \right),
\end{aligned}
\end{equation}
and their corresponding spatial derivatives are:
\begin{equation}
\begin{aligned}
    q_{1,x}^l &= \frac{1}{h} \left( -3Q_{j-1} + 3Q_j - 2Q_{x,j-1}h \right), \\
    q_{2,x}^l &= \frac{1}{h} \left( Q_{j+1} - Q_j \right), \\
    q_{3,x}^l &= \frac{1}{h} \left( Q_{j+1} - Q_j \right).
\end{aligned}
\end{equation}
The smoothness indicators \(IS_{k}^l\) for the left-side interface reconstruction are defined as:
\begin{equation}
\begin{aligned}
    IS_{1}^l &= \frac{13}{3} \left( Q_j - Q_{j-1} - Q_{x,j-1}h \right)^2 + \left( 2(Q_j - Q_{j-1}) - Q_{x,j-1}h \right)^2, \\
    IS_{2}^l &= \frac{13}{12} \left( Q_{j-1} - 2Q_j + Q_{j+1} \right)^2 + \frac{1}{4} \left( Q_{j+1} - Q_{j-1} \right)^2, \\
    IS_{3}^l &= \frac{13}{12} \left( Q_j - 2Q_{j+1} + Q_{j+2} \right)^2 + \frac{1}{4} \left( 3Q_j - 4Q_{j+1} + Q_{j+2} \right)^2.
\end{aligned}
\end{equation}

\subsection{Reconstruction at interfaces (right side)}
\label{appendix:C2}

For the right-side reconstruction at the grid interface \(j-1/2\), the sub-stencil polynomials are:
\begin{equation}
\begin{aligned}
    q_{1}^r &= \frac{1}{4} \left( -5Q_{j+2} + 9Q_{j+1} + 3Q_{x,j+2}h \right), \\
    q_{2}^r &= \frac{1}{8} \left( -Q_{j+2} + 6Q_{j+1} + 3Q_j \right), \\
    q_{3}^r &= \frac{1}{8} \left( -Q_{j-1} + 6Q_j + 3Q_{j+1} \right),
\end{aligned}
\end{equation}
and their spatial derivatives are:
\begin{equation}
\begin{aligned}
    q_{1,x}^r &= \frac{1}{h} \left( -3Q_{j+1} + 3Q_{j+2} - 2Q_{x,j+2}h \right), \\
    q_{2,x}^r &= \frac{1}{h} \left( Q_{j+1} - Q_j \right), \\
    q_{3,x}^r &= \frac{1}{h} \left( Q_{j+1} - Q_j \right).
\end{aligned}
\end{equation}
The corresponding smoothness indicators \(IS_{k}^r\) are:
\begin{equation}
\begin{aligned}
    IS_{1}^r &= \frac{13}{3} \left( Q_{j+1} - Q_{j+2} + Q_{x,j+2}h \right)^2 + \left( 2(Q_{j+2} - Q_{j+1}) - Q_{x,j+2}h \right)^2, \\
    IS_{2}^r &= \frac{13}{12} \left( Q_j - 2Q_{j+1} + Q_{j+2} \right)^2 + \frac{1}{4} \left( Q_{j+2} - Q_j \right)^2, \\
    IS_{3}^r &= \frac{13}{12} \left( Q_{j-1} - 2Q_j + Q_{j+1} \right)^2 + \frac{1}{4} \left( 3Q_{j+1} - 4Q_j + Q_{j-1} \right)^2.
\end{aligned}
\end{equation}

\subsection{Reconstruction at nodes}
\label{appendix:C3}

For the derivative reconstruction at the solution node \(j\), the sub-stencil polynomials \(q_{k}^N\) are:
\begin{equation}
\begin{aligned}
    q_{1,x}^N &= \frac{1}{h} \left( -2Q_{j-1} + 2Q_j - Q_{x,j-1}h \right), \\
    q_{2,x}^N &= \frac{1}{h} \left( -2Q_j + 2Q_{j+1} - Q_{x,j+1}h \right), \\
    q_{3,x}^N &= \frac{1}{2h} \left( Q_{j+1} - Q_{j-1} \right).
\end{aligned}
\end{equation}
The smoothness indicators \(IS_{k}^N\) for the node reconstruction are evaluated as:
\begin{equation}
\begin{aligned}
    IS_{1}^N &= \frac{13}{3} \left( Q_j - Q_{j-1} - Q_{x,j-1}h \right)^2 + \left( 2(Q_j - Q_{j-1}) - Q_{x,j-1}h \right)^2, \\
    IS_{2}^N &= \frac{13}{3} \left( Q_{j+1} - Q_j - Q_{x,j+1}h \right)^2 + \left( 2(Q_{j+1} - Q_j) - Q_{x,j+1}h \right)^2, \\
    IS_{3}^N &= \frac{13}{12} \left( Q_{j-1} - 2Q_j + Q_{j+1} \right)^2 + \frac{1}{4} \left( Q_{j+1} - Q_{j-1} \right)^2.
\end{aligned}
\end{equation}

\section{FORTRAN program of GENO method}
\label{appendix:D}

The following Fortran subroutine implements the GENO nonlinear reconstruction.

\begin{verbatim}
subroutine Func_Geno_reconstruction(IS, p, Wr)
    !-----------------------------------------------------------------
    ! GENO nonlinear reconstruction (value only).
    !
    ! Inputs:
    !   IS(1:3) -- smoothness indicators of the three sub-stencils
    !   p(1:4)  -- candidate polynomial values at the interface:
    !              p(1:3): low-order sub-stencil polynomials
    !              p(4)  : high-order polynomial
    ! Output:
    !   Wr      -- reconstructed value
    !-----------------------------------------------------------------
    implicit none

    real(8), intent(in)  :: IS(3), p(4)
    real(8), intent(out) :: Wr

    real(8) :: tau, w(2), wh, a(3)
    real(8), parameter :: eps = 1.0d-15, eps_a = 1.0d-6, s = 20.0d0
    real(8), parameter :: d(3) = (/ 1.0d0, 8.0d0, 1.0d0 /)

    ! Step 1: path-function-based hybrid weight
    tau  = 0.5d0 * abs(IS(2) + IS(3) - 2.0d0*IS(1))
    w(1) = 1.0d0 + (tau / (maxval(IS) + eps))**3
    w(2) = 1.0d0 + (tau / (minval(IS) + eps))**3
    w    = 2.0d0 * w / (w(1) + w(2))
    wh   = tanh(s * w(1)) / tanh(s)          ! geno transition

    ! Step 2: GENO nonlinear weights
    a = d / (IS + eps_a)**2
    a = a / sum(a)

    ! Step 3: convex combination of low- and high-order polynomials
    Wr = (1.0d0 - wh) * sum(a * p(1:3)) + wh * p(4)

end subroutine Func_Geno_reconstruction
\end{verbatim}

\section*{References}
\bibliographystyle{ieeetr}
\bibliography{compact-gks}

\end{document}